\documentclass[letterpaper,journal,onecolumn,oneside,final,12pt]{IEEEtran}
\usepackage{amssymb} 
\usepackage{color}
\definecolor{dkgreen}{rgb}{0,.5,0}
\definecolor{dkred}{rgb}{.5,0,0}
\definecolor{dkcyan}{rgb}{0,.5,.5}
\definecolor{dkmgnt}{rgb}{.5,0,.5}
\usepackage[dvips,colorlinks,linkcolor=dkred,pagecolor=red,filecolor=dkcyan,citecolor=dkgreen,urlcolor=dkmgnt%
,bookmarksnumbered=true%
,pdftitle=On\ $Z_{2^k}$-Dual\ Binary\ Codes%
,pdfsubject=Coding\ Theory%
,pdfauthor=Denis\ S.\ Krotov%
,pdfkeywords=Gray\ map;\ Hadamard\ codes;\ MacWilliams\ identity;\ perfect\ codes;\ $Z_{2^k}$-linearity%
,hypertexnames=false%
]{hyperref} 

\providecommand\hyperref[2][]{#2} 
\providecommand\triangleq{{\stackrel{\scriptscriptstyle\triangle}{{}={}}}} 

\renewcommand{\QED}{\QEDclosed} 
\newcommand{\proofr}{\begin{proof}} 
\newcommand{\proofrem}[1]{\par\indent\indent \emph{#1:}} 
\newcommand{\proofend}{\end{proof}}  
\newcommand{\proofremend}{\hfill\QED}  

\newtheorem{theorem}{Theorem}[section]
\newtheorem{lemma}[theorem]{Lemma}
\newtheorem{pro}[theorem]{Proposition}
\newtheorem{corol}[theorem]{Corollary}
\newtheorem{note}[theorem]{Remark}
\newtheorem{examp}[theorem]{Example}

\newcommand{\bpro}{\begin{pro}}
\newcommand{\bprorem}[1]{\begin{pro}[#1]}
\newcommand{\epro}{\end{pro}}

\newcommand{\bcor}{\begin{corol}}
\newcommand{\ecor}{\end{corol}}

\newcommand{\bnote}{\begin{note}}
\newcommand{\enote}{\end{note}}

\newcommand{\btheorem}{\begin{theorem}}
\newcommand{\etheorem}{\end{theorem}}

\newcommand{\blemma}{\begin{lemma}}
\newcommand{\blemmarem}[1]{\begin{lemma}[#1]}
\newcommand{\elemma}{\end{lemma}}

\newcommand{\bexam}{\begin{examp}}
\newcommand{\eexam}{\end{examp}}

\newcommand\1{$\overline 1$}

\begin{document}

\title{On $Z_{2^k}$-Dual Binary Codes%
\thanks{This is author's version of the correspondence in the
IEEE Transactions on Information Theory 53(4) 2007, 1532--1537, Digital
Object Identifier \href{http://dx.doi.org/10.1109/TIT.2007.892787}{10.1109/TIT.2007.892787}, \copyright
2007 IEEE.}
\thanks{Some results of this paper were presented at
the 4th International Workshop on Optimal Codes and Related Topics
OC\,2005 (Pamporovo, Bulgaria, June 2005).}%
}
\renewcommand\footnotemark{ }
\author{Denis~S.~Krotov%
\thanks{D.\,S.\,Krotov is with the Sobolev Institute of Mathematics, Novosibirsk, Russia
(e-mail: {\tt krotov@math.nsc.ru}).}
}

\markboth
{Krotov: On $Z_{2^k}$-Dual Binary Codes}
{Krotov: On $Z_{2^k}$-Dual Binary Codes}

\maketitle

\begin{abstract}
A new generalization of the Gray map is introduced.
The new generalization $\Phi: Z_{2^k}^n \to Z_{2}^{2^{k-1}n}$
is connected with the known generalized Gray map $\varphi$
in the following way:
if we take two dual linear $Z_{2^k}$-codes and construct
binary codes from them using the generalizations $\varphi$ and $\Phi$ of the Gray map,
then the weight enumerators of the binary codes obtained will satisfy
the MacWilliams identity.
The classes of
$Z_{2^k}$-linear Hadamard codes and co-$Z_{2^k}$-linear
extended 1-perfect codes are described,
where co-$Z_{2^k}$-linearity means that the code can be obtained from
a linear $Z_{2^k}$-code with the help of the new generalized Gray map.
\end{abstract}

\section{Introduction}\label{S_Intro}
{A}{s} discovered in \cite{Nechaev82,Nechaev:Kerdock_cyclic_form}
certain nonlinear binary codes can be represented as linear codes over $Z_4$.
The variant of this representation founded in \cite{HammonsOth:Z4_linearity}
use the mapping $\phi: 0\to 00,\ 1\to 01,\ 2\to 11, 3\to 10$,
which is called the {\it Gray map},
to construct binary so-called $Z_4$-linear codes from linear quaternary codes.
The main property of $\phi$ from this point of view is that it is an isometry
between $Z_4$ with the Lee metric and $Z_2^2$ with the Hamming metric.
In \cite{Carlet:Z2k-linear} (and in \cite{HonNech} in more general form)
the Gray map is generalized to construct $Z_{2^k}$-linear codes.
The generalized Gray map
(say $\varphi$; see \hyperref[S_Z-lin0]{Subsection~}\ref{S_Z-lin0} for recalling basic facts
on the generalized Gray map)
is an isometric imbedding of $Z_{2^k}$ with the metric
specified by the homogeneous weight \cite{ConstHeise}
into $Z_2^{2^{k-1}}$ with the Hamming metric.

In this paper we introduce another generalization $\Phi$ of the Gray map
(\hyperref[S_co-Z-lin]{Subsection~}\ref{S_co-Z-lin}).
This generalization turns out to be dual to the previous in the following sense.
If $\cal C$ and ${\cal C}^\bot$ are dual linear $Z_{2^k}$-codes,
then the binary $Z_{2^k}$-linear
code $\varphi({\cal C})$ and the {\it co-$Z_{2^k}$-linear}
code $\Phi({\cal C}^\bot)$ are formally dual.
The formal duality is that the weight enumerators
of these two codes satisfy the MacWilliams identity (\hyperref[S_dual]{Section~}\ref{S_dual});
note that these codes, in general, can be nonlinear, in which case they cannot be
dual in the usual sense, as subspaces of the binary vector space.
So, we solve the problem of duality
for $Z_{2^k}$-linear binary codes: to relate the weight enumerators of the images
of dual linear codes over $Z_{2^k}$. This problem
cannot be solved using only the standard generalized Gray map $\varphi$,
because, as noted in \cite{Carlet:Z2k-linear},
the weight enumerators of $\varphi({\cal C})$ and $\varphi({\cal C}^\bot)$
are in general not related, contrarily to the case of $Z_4$-linear codes.

In \cite{HonNech}, it is shown that binary (and, in general, nonbinary)
codes can be represented as group codes over different groups.
Such representations use a special \emph{scaled isometry}
($\varphi$ is a partial case of such isometry),
which acts isometrically from
some module with a specially defined metric to the binary (or, in general, nonbinary)
Hamming space.
With such approach, each module element corresponds to one codeword;
$Z_{2^k}$-linear codes are a partial case of that approach.

In our approach, every module element (word) corresponds to some set of codewords,
which is a Cartesian product of the sets corresponding,
by a special mapping, to the symbols of the module word.
In the constructive part,
this approach can be represented
by the generalized concatenation construction \cite{Zin1976:GCC},
but the resulting code distance differs
from the constructive distance of generalized concatenated codes.
Calculating distance, we use some isometrical properties of the mapping, which
can be considered as an analog of a scaled isometry.
The map $\Phi$ considered in this paper allows to construct only codes with distance maximum $4$,
but in general the approach may have a larger potential.

The new way to generalize the Gray map (\hyperref[S_co-Z-lin]{Subsection~}\ref{S_co-Z-lin})
and its duality (\hyperref[S_dual]{Section~}\ref{S_dual}) to the ``old way''
are the main results of the first part of the paper.
The second part (\hyperref[S_perf]{Sections~}\ref{S_perf}\hyperref[S_Had]{--}\ref{S_nonexist})
is devoted to constructions of
co-$Z_{2^k}$-linear and $Z_{2^k}$-linear codes with parameters of
extended $1$-perfect and Hadamard codes.
There is a series
(e.\,g.
\cite{RifPuj:1997},
\cite{BorRif:1999},
\cite{Kro:2000:Z4_Perf},
\cite{Kro:2001:Z4_Had_Perf},
\cite{PheRif:2002},
\cite{BorPheRif:2003},
\cite{PheRifVil:2004ACCT:additive_hadamard},
\cite{PheRifVil:Hadamard-rank-kernel},
\cite{GuptaOth2005})
of papers that concern $1$-perfect (or extended $1$-perfect)
and Hadamard codes which are linear in some nonstandard sense
(note that for each considered values of parameters
only exists one linear code, up to equivalence).
Using the two generalized $\varphi$, $\Phi$ Gray maps, we construct a wide class of such codes.
Some natural questions on the constructed codes remains open for future researching:
which of these codes are equivalent to each other or to other known codes;
which of them are propelinear (see the definition in \cite{RifPuj:1997}) or
at least transitive
(see \cite{Sol:2005transitive} for the definition
and some constructions of transitive $1$-perfect codes);
to establish the bounds on the dimension of the kernel and rank
for these codes
(the dimension of the kernel and rank are good measures of linearity
of nonlinear codes in the classical binary sense.
For $Z_4$-linear extended $1$-per\-fect and Hadamard codes these parameters were calculated in
\cite{Kro:2000:Z4_Perf}, \cite{Kro:2001:Z4_Had_Perf}, \cite{BorPheRif:2003}, \cite{PheRifVil:2004ACCT:additive_hadamard}, \cite{PheRifVil:Hadamard-rank-kernel});
and so on.

In \hyperref[S_perf]{Section~}\ref{S_perf} we construct
co-$Z_{2^k}$-linear extended $1$-perfect $(n,2^n/2n,4)$ codes.
As noted above, the code distance of a co-$Z_{2^k}$-linear code
cannot exceed $4$ if $k>2$ (see \hyperref[DcoLin]{Lemma~}\ref{DcoLin}).
So, extended $1$-perfect codes are the best examples.

As a corollary, the dual distance of a $Z_{2^k}$-linear code ($k>2$) also
cannot exceed $4$. Best examples are codes with large code distance.
In \hyperref[S_Had]{Section~}\ref{S_Had} we construct $Z_{2^k}$-linear codes with
parameters $(n,2n,n/2)$, i.\,e., Hadamard codes.
They are
$Z_{2^k}$-dual to the extended $1$-perfect codes
from \hyperref[S_perf]{Section~}\ref{S_perf};
i.\,e., the $\varphi$ and, respectively, $\Phi$ preimage of the constructed
extended $1$-perfect and Hadamard codes are dual $Z_{2^k}$ codes.

The introduced constructions of codes are a generalization of
the constructions of $Z_4$-linear extended $1$-perfect and Hadamard codes
\cite{Kro:2000:Z4_Perf,Kro:2001:Z4_Had_Perf} (see also \cite{BorPheRif:2003}).

It is natural to ask whether our construction of
co-$Z_{2^k}$-linear extended $1$-perfect
and $Z_{2^k}$-linear Hadamard codes is complete or not.
In \hyperref[S_nonexist]{Section~}\ref{S_nonexist} we
show that the answer is yes;
i.\,e., all co-$Z_{2^k}$-linear extended $1$-perfect
codes and all $Z_{2^k}$-linear Hadamard codes are equivalent to codes constructed
in \hyperref[S_perf]{Sections~}\ref{S_perf} \hyperref[S_Had]{and~}\ref{S_Had}.

\section[Definitions And Basic Facts]{Definitions And Basic Facts. Two Generalizations Of Gray Map}\label{S_def}
{W}{e} will use the following common notations.
An \emph{$(n,M,d)$ binary code} is a cardinality $M$ subset of
$Z_2^n = {Z_2\times\ldots\times Z_2}$
with the distance at least
$d$ between every two different elements (\emph{codewords});
$n$ is called the code \emph{length} and $d$ is called the \emph{code distance}.
An $(n,2n,n/2)$ binary code is called an \emph{Hadamard code}; such codes exist
if and only if Hadamard $n\times n$ matrices exist, i.\,e.,
at least $n\equiv 0\mbox{ mod }4$ if $n>2$, see e.\,g. \cite{MWS}.
A linear Hadamard code is known as
a first order Reed-Muller code;
$n$ is a power of two is this case.
An $(n,2^n/2n,4)$ binary code is called an \emph{extended $1$-perfect code};
such codes exist if and only if $n$ is a power of two.
Such a code is always an extension of a $(n-1,2^n/2n,3)$ code, which is called
$1$-per\-fect; so, studying extended $1$-perfect codes is an alternative way
to study $1$-perfect codes. For each admissible length there is only one,
up to coordinate permutation, linear (extended) $1$-perfect code and,
if $n>8$, many nonlinear ones;
the classification of all such codes is currently  an open problem.
In the following two subsections
we will define the basic concepts, which will be used throughout the paper:
two weight functions on $Z_{2^k}$, two corresponding metrics,
two generalized Gray maps $\varphi$ and $\Phi$,
the concepts of $Z_{2^k}$-li\-ne\-ar and co-$Z_{2^k}$-li\-ne\-ar codes;
we will formulate simple but fundamental claims on isometric properties of $\varphi$ and $\Phi$.
\subsection[$Z_{2m}$-Linear Codes]{$Z_{2m}$-Linear Codes}\label{S_Z-lin0}

This subsection recalls basic definitions and facts about the
generalized Gray map and $Z_{2^k}$-linear codes \cite{Carlet:Z2k-linear,HonNech}.

Assume $m \geq 2$ is an integer such that there exist an Hadamard
$m\times m$ matrix.
Let $A\subset Z_2^{m}$ be a Hadamard $(m,2m,m/2)$ code.
Assume that $A=\{a_0,a_1, \ldots ,a_{2m-1}\}$, where $a_0$
is the all-zero codeword
and $a_i+a_{i+m}$ is the all-one codeword for each $i$ from $0$ to $m-1$.
Define the {\em generalized Gray map} $\varphi : Z_{2m}^n \to Z_{2}^{mn}$ by the rule
$$
\varphi(x_1,\ldots,x_n)\triangleq (a_{x_1},\ldots,a_{x_n}).
$$
Let the weight function $wt^*:Z_{2m}\to R^+$ be given by
$$
wt^*(x) \triangleq
\cases{
0& \mbox{if $x=0$}, \cr
m& \mbox{if $x=m$}, \cr
m/2& \mbox{otherwise}.}
$$
Whenever $m$ is a power of two, the weight $wt^*$ is the \emph{homogeneous weight} introduced in \cite{ConstHeise}
for more general class of rings. The corresponding distance $d^*$ on $Z_{2m}^n$
is defined by the standard way:
$$ d^*((x_1,\ldots,x_n),(y_1,\ldots,y_n))
\triangleq
\sum_{i=1}^n wt^*(y_i-x_i).
$$

\bprorem{\cite{HonNech}} \label{P_iso}
The mapping $\varphi$ is an isometric embedding of $(Z_{2m}^n,d^*)$ into
$(Z_2^{mn},d)$, i.\,e.,
$d^*(\bar x,\bar y)=d(\varphi(\bar x),\varphi(\bar y))$, where $d(\cdot,\cdot)$
is the Hamming distance.
\epro

If ${\cal C}\subseteq Z_{2m}^n$, then
$$
\varphi({\cal C})\triangleq \{ \varphi(\bar x) \,|\, \bar x\in {\cal C} \}.
$$

We will say that ${\cal C}$ is a $(n,M,d)^*$ code if 
${\cal C}\subseteq Z_{2m}^n$,
$|{\cal C}|=M$, and the $d^*$-distance between
any two different elements of ${\cal C}$ is not less than $d$.

\bcor\label{Cor1}
If ${\cal C}$ is a $(n,M,d)^*$ code, then $\varphi({\cal C})$
is a binary $(mn,M,d)$ code.
\ecor

A binary code is called {\em $Z_{2m}$-linear} if 
its coordinates can be arranged in such a way that
it is the image of a linear $Z_{2m}$-code by $\varphi$.
As noted in \cite{Carlet:Z2k-linear}, the length of a $Z_{2m}$-linear code
must be a multiple of $m$ and all the weights of its codewords
must be multiples of $m/2$.

\bnote
The map $\varphi$ and the concept
of $Z_{2m}$-linearity given above depends on the choice of
a Hadamard code $A$ and the only restriction on $m$ is that exists
a Hadamard $m\times m$ matrix.
In fact, originally {\rm \cite{Carlet:Z2k-linear}}
(and in the binary subcase of \cite[Section D]{HonNech})
$m=2^{k-1}$ and the code $A$ is fixed as
$R(1,k-1)$, the Reed-Muller code of order $1$.
\enote

\subsection[Co-$Z_{2^k}$-Linear Codes]{Co-$Z_{2^k}$-Linear Codes}
\label{S_co-Z-lin}
In this subsection we will introduce another approach
to construct binary codes from linear codes over $Z_{2^k}$.

Firstly, we will introduce a new way to generalize the Gray map, defining a map
$\Phi : Z_{2m}^n \to 2^{Z_{2}^{mn}}$,
where $2^{Z_{2}^{mn}}$ denotes the set of the subsets of $Z_{2}^{mn}$.
Then we will define a weight function $wt^\diamond:Z_{2m}\to R^+$
and the corresponding distance
$d^\diamond$ and we will show (\hyperref[DcoLin]{Lemma~}\ref{DcoLin})
relations between the $d^\diamond$-distance of a code ${\cal C}\subset Z_{2m}^n$
and the Hamming distance of the code $\Phi({\cal C})$.
Finally, we will introduce the concept of co-$Z_{2^k}$-linear codes.

Put $m=2^{k-1}$. Let $\{ H_0,\ldots,H_{2m-1} \}$ be a partition of $Z_2^{m}$ into
extended $1$-perfect $(m,2^m/2m,4)$ codes
(for example, we can take $H_0$ as the extended Hamming code and $H_0,\ldots,H_{2m-1}$ as its cosets).
Moreover, we assume that $H_0$ contains the all-zero word $\bar 0$ and
$H_j$ is an even (odd) weighted if and only if $j$ is even (odd).

Define the map $\Phi : Z_{2m}^n \to 2^{Z_{2}^{mn}}$ by the rule
$$
\Phi(x_1,\ldots,x_n)\triangleq H_{x_1}\times\ldots\times H_{x_n}.
$$
If ${\cal C}\subseteq Z_{2m}^n$, then\\[-10pt]
$$
\Phi({\cal C})\triangleq \bigcup_{\bar x\in {\cal C}} \Phi(\bar x).
$$

\bexam
Let $k=3$, $H_0=\{0000,1111\}$, \ldots,
$H_2=\{1100,0011\}$, \ldots,
 $H_7=\{0001,1110\}$. Then
\begin{IEEEeqnarray*}{rCrCl}
\Phi(2\,0\,7)&=&H_2\times H_0\times H_7
&=&\{\, 1100\,0000\,0001,\ 1100\,0000\,1110,\ 1100\,1111\,0001,\ 1100\,1111\,1110,\\
&&&&\phantom{\{\,} 0011\,0000\,0001,\ 0011\,0000\,1110,\ 0011\,1111\,0001,\ 0011\,1111\,1110\,\}.
\end{IEEEeqnarray*}
\eexam
\hspace{1.0ex}

Define the weight function $wt^\diamond:Z_{2m}\to R^+$ as
$$
wt^\diamond(x)\triangleq
\cases{
0&\mbox{if $x=0$},\cr
1&\mbox{if $x$ is odd},\cr
2&\mbox{if $x\neq 0$ is even}.}
$$
The corresponding distance $d^\diamond$ on $Z_{2m}^n$ is defined by the standard way:
$$ d^\diamond((x_1,\ldots,x_n),(y_1,\ldots,y_n))
\triangleq
\sum_{i=1}^n wt^\diamond(y_i-x_i).
$$

We will say that ${\cal C}$ is an $(n,M,d)^\diamond$ code if 
${\cal C}\subseteq Z_{2m}^n$,
$|{\cal C}|=M$, and the $d^\diamond$-distance between
any two different elements of ${\cal C}$ is not less than $d$.

\blemma\label{DcoLin}
Let $m\geq 4$.
If ${\cal C}\subseteq Z_{2m}^n$ is an $(n,M,d)^\diamond$ code,\\
then $\Phi({\cal C})$ is a binary $(mn,M (\frac{2^m}{2m})^n,\min(4,d))$ code.
\elemma
The proof is straightforward, and we omit it.
%
%
%
%
%
%

We call a binary code {\it co-$Z_{2m}$-linear} if 
its coordinates can be arranged so that
it is the image of a linear $Z_{2m}$-code by the map $\Phi$.

\bnote
Co-$Z_{2m}$-linear codes can be considered as a partial case
of generalized concatenated codes \cite{Zin1976:GCC}.
For this special case the code distance given by \hyperref[DcoLin]{Lemma~}\ref{DcoLin}
may be better than the code distance guaranteed by the
generalized concatenation construction.
\enote

\section[$Z_{2^k}$-Duality Of Binary Codes]{%
$Z_{2^k}$-Duality Of Binary Codes}
\label{S_dual}
{I}{n} this section we will show (\hyperref[duality]{Theorem~}\ref{duality})
that if two binary codes are obtained from
dual linear $Z_{2^k}$-codes by, both, $\varphi$ and $\Phi$
generalizations of Gray map,
then these codes are formally dual, i.\,e., their weight enumerators satisfy the
MacWilliams identity.

Let $C$ be a $Z_{2^k}$-linear code, the image by $\varphi$
of a linear $Z_{2^k}$-code $\cal C$ of length $n$.
Let ${\cal C}^\bot$ be the linear $Z_{2^k}$-code dual to ${\cal C}$.
Let $SW_{{\cal C}^\bot}(X,Z,T)$ be the polynomial obtained from the complete
weight enumerator $W_{{\cal C}^\bot}(X_0,X_1,\ldots,X_{2^k-1})$ of ${\cal C}^\bot$
by identifying to $Z$ (respectively to $T$)
all the $X_j$'s such that $j$ is odd (respectively, $j$ is even $\neq 0$).

\blemmarem{\cite{Carlet:Z2k-linear}}\label{CC}
Then we have
\begin{IEEEeqnarray*}{rRl}
W_C(X,Y)&
=&\frac{1}{|{\cal C}^\bot|}SW_{{\cal C}^\bot} \big(\
X^{2^{k-1}}\hspace{-0.8ex}+Y^{2^{k-1}}\hspace{-0.8ex}+(2^k-2)(X Y)^{2^{k-2}}\hspace{-0.4ex},\ \nonumber \\
&&\phantom{\frac{1}{|{\cal C}^\bot|}SW_{{\cal C}^\bot} \big(\ }
X^{2^{k-1}}\hspace{-0.8ex}-Y^{2^{k-1}},\  \nonumber \\
&&\phantom{\frac{1}{|{\cal C}^\bot|}SW_{{\cal C}^\bot} \big(\ }
X^{2^{k-1}}\hspace{-0.8ex}+Y^{2^{k-1}}\hspace{-0.8ex}-2(X Y)^{2^{k-2}}\ \big).\nonumber \\
\end{IEEEeqnarray*}
\elemma

The following lemma is a known fact about the weight distribution of
extended $1$-perfect binary codes.
\blemma\label{W_H}
Let $H$ be an extended $1$-perfect $(m,2^m/2m,4)$ code;
then \\
{\rm a)}
\(
\frac{1}{|H|}
W_H(X+Y,X-Y)=X^{m}+Y^{m}+(2m-2)(X Y)^{m/2}
\) \ if $\bar 0\in H$,
\\[0pt]
{\rm b)}
\(
\frac{1}{|H|}
W_H(X+Y,X-Y)=X^{m}-Y^{m}
\) \ if $H$ is odd-weight,
\\[0pt]
{\rm c)}
\(
\frac{1}{|H|}
W_H(X+Y,X-Y)=X^{m}+Y^{m}-2(X Y)^{m/2}
\) \ if $H$ is an even-weight code and $\bar 0\not\in H$.%
\elemma

\bpro\label{W-SW}
Let $\cal C$ be a $Z_{2^k}$-code and
$\tilde C=\Phi({\cal C})$; then
$$
W_{\tilde C}(X,Y)=SW_{\cal C}(W_{H_0}(X,Y), W_{H_1}(X,Y), W_{H_2}(X,Y)),
$$
where $\{ H_0,\ldots,H_{2m-1} \}$ is the partition of $Z_2^{m}$ into
extended $1$-perfect codes from the definition of $\Phi$,
$0\in H_0$, $H_1$ is an odd-weight code,
$H_2$ is an even-weight code, $\bar 0\not\in H_2$
\epro
\proofr
The statement follows almost immediately from the definition of the map $\Phi$.
Indeed, each codeword $\bar z=(z_1,\ldots,z_n)\in {\cal C}$ adds
$SW_{z_1}\cdot\ldots\cdot SW_{z_n}$ to $SW_{\cal C}(X,Z,T)$, where
\begin{IEEEeqnarray*}{rCl"l}
  SW_{0} &\triangleq& X, \\
  SW_{2j+1}  &\triangleq& Z,& j=0, \ldots ,m-1,\\
  SW_{2j}  &\triangleq& T, & j=1, \ldots ,m-1.
\end{IEEEeqnarray*}
On the other hand, as follows from the definition of the map $\Phi$,
each codeword $\bar z=(z_1,\ldots,z_n)\in {\cal C}$ adds
$W_{H_{z_1}}(X,Y)\cdot\ldots\cdot W_{H_{z_n}}(X,Y)$ to $W_{\tilde C}(X,Y)$.
The relations
\begin{IEEEeqnarray*}{rCl"l}
  W_{H_{2j+1}}(X,Y)&=&W_{H_{1}}(X,Y), \qquad j=0, \ldots ,m-1,\\
  W_{H_{2j}}(X,Y)&=&W_{H_{2}}(X,Y), \qquad j=1, \ldots ,m-1
\end{IEEEeqnarray*}
follow from \hyperref[W_H]{Lemma~}\ref{W_H} and conclude the proof.
\proofend

The following theorem is the main result of this section.
\btheorem\label{duality}
Let ${\cal C}$ and ${\cal C}^\bot$ be dual linear $Z_{2^k}$-codes,
$C=\varphi({\cal C})$ and $\tilde C_\bot=\Phi({\cal C}^\bot)$.
Then the codes $C$ and $\tilde C_\bot$ are formally dual, i.\,e.,
\(
 W_C(X,Y)=\frac{1}{|\tilde C_\bot|}W_{\tilde C_\bot}(X+Y,X-Y).
\)
\etheorem
\proofr
By \hyperref[CC]{Lemmas~}\ref{CC} \hyperref[W_H]{and~}\ref{W_H}
\begin{IEEEeqnarray*}{rCl}
W_C(X,Y)
&=&\frac{1}{|{\cal C}^\bot|}SW_{{\cal C}^\bot} \bigg(
\frac{1}{|H_0|}W_{H_0}(X{+}Y,X{-}Y),
\frac{1}{|H_1|}W_{H_1}(X{+}Y,X{-}Y),
\frac{1}{|H_2|}W_{H_2}(X{+}Y,X{-}Y)
\bigg)\\
&=&\frac{1}{|{\cal C}^\bot|}\left(\frac{2m}{2^m}\right)^{\!\!n}SW_{{\cal C}^\bot} \Big(
W_{H_0}(X{+}Y,X{-}Y),
W_{H_1}(X{+}Y,X{-}Y),
W_{H_2}(X{+}Y,X{-}Y)
\Big).
\end{IEEEeqnarray*}
It remains to note that by \hyperref[DcoLin]{Lemma~}\ref{DcoLin} we have
${|{\cal C}^\bot|}\left(\frac{2^m}{2m}\right)^n =|\tilde C_\bot|$
and by \hyperref[W-SW]{Proposition~}\ref{W-SW} we obtain:
\begin{IEEEeqnarray*}{rlCl}
SW_{{\cal C}^\bot} \big(W_{H_0}(X{+}Y,X{-}Y),
W_{H_1}(X{+}Y,X{-}Y),
W_{H_2}(X{+}Y,X{-}Y){\big)}
{}={}&W_{\tilde C_\bot}(X{+}Y,X{-}Y).
\vspace{-4.1ex} 
\end{IEEEeqnarray*}
\vspace{2ex} 
\proofend
In \hyperref[S_perf]{Sections~}\ref{S_perf} \hyperref[S_Had]{and~}\ref{S_Had}
we will construct two $Z_{2^k}$-dual classes of codes:
co-$Z_{2^k}$-linear extended $1$-perfect codes
and $Z_{2^k}$-linear Hadamard codes.
In the \hyperref[S_nonexist]{last section}
we will show the completeness of the constructions.
\section[Co-$Z_{2^k}$-Linear Extended $1$-Perfect Codes]{%
Co-$Z_{2^k}$-Linear Extended Perfect Codes}
\label{S_perf}

{T}{his} section concerns extended $1$-perfect codes.
We first introduce the concept of {\em \1-perfect code}, which is a generalization of
the concept of extended $1$-perfect binary code to some nonbinary cases.
As in the case of $1$-per\-fect codes, the existence of \1-per\-fect codes
in different spaces is independently interesting.
Then, in \hyperref[ss_b]{Subsection~}\ref{ss_b},
we  construct a class of {\em \1-perfect codes} in $(Z_{2^k}^n,d^\diamond)$.
In \hyperref[ss_c]{Subsection~}\ref{ss_c} we  summarize:
the images of such codes under $\Phi$ are co-$Z_{2^k}$-linear extended $1$-perfect codes.
In \hyperref[ss_d]{Subsection~}\ref{ss_d}
we give examples of \1-per\-fect codes in $Z_{2m}$ where $m$ is not a power of two.

\subsection[$1'$-Perfect Codes]{\1-Perfect Codes} \label{ss_a}
Let $G=(V,E)$ be a regular bipartite graph with parts $V_{ev}$, $V_{od}$.
A subset $C\subseteq V_{ev}$ is called a {\em \1-perfect code} if 
for each $\bar x\in V_{od}$ there exist exactly one $\bar c\in C$
adjacent with $\bar x$.
It is not difficult to see that a \1-perfect code is an optimal distance $4$ code,
i.\,e., its cardinality is maximum among the distance $4$ codes in the same space.
In particular, such a code is a $3$-diameter perfect code in the sense \cite{AhlAydKha}.
\1-perfect codes in a binary Hamming space are known as extended $1$-perfect codes.
The following test can be considered as an alternative definition of a \1-perfect code.
\bpro\label{1perf}
Assume $G=(V,E)$ is a regular bipartite graph  of degree $r>0$
and $d^{G}$ is the graph distance in $V$.
Then $C\subset V$ is a \1-perfect code if and only if
$|C|=|V|/2r$ and $d^{G}(\bar c_1,\bar c_2)\geq 4$
for each different $\bar c_1, \bar c_2\in C$.
\epro
\proofr
{\em Only if:} Assume $C\subset V$ is a \1-perfect code.
Then by the definition it is a subset of a part $V_{ev}$
of the graph $(V,E)$ and the graph distance between any two elements of $C$ is even.

On the other hand, if $d^{G}(\bar c_1, \bar c_2)=2$ for some $\bar c_1, \bar c_2\in C$,
then there exist an element $\bar x\in V$ such that
$d^{G}(\bar x, \bar c_1)=1$ and $d^{G}(\bar x, \bar c_2)=1$.
The existence of such element contradicts to the definition of \1-perfect code.
Consequently, $d^{G}(\bar c_1, \bar c_2)\geq 4$ for each different $\bar c_1, \bar c_2\in C$.

Each element of $V_{od}\triangleq V\backslash V_{ev}$ is adjacent with exactly one element of $C$.
On the other hand, each element of $C$ is adjacent with exactly $r$ elements of $V_{od}$.
Therefore, $|C|=|V_{od}|/r=|V|/2r$.

{\em If:}
Assume $C$ is a code of cardinality $|V|/2r$ and distance at least $4$.
Denote
$$
C_1 \triangleq \{\bar x\,|\,d^{G}(\bar x,C)=1\}.
$$
Each element of $C_1$ is adjacent with exactly one element of $C$
(otherwise the code distance of $C$ is not more than $2$).
Consequently, $|C_1|=r|C|=|V|/2$.

(*) {\em We claim that $C_1$ does not contain two adjacent elements.}
Assume the contrary.
Then the graph $G$ contains a chain $(\bar c, \bar x, \bar y, \bar c')$,
where $\bar c, \bar c'\in C$, $\bar x, \bar y \in C_1$.
The case $\bar c=\bar c'$ contradicts to the bipartiteness of the graph $G$.
The case $\bar c\neq\bar c'$ contradicts to the code distance of $C$.
The claim (*) is proved.

Since $G$ is a regular graph, the other half of vertices $V \backslash C_1$
also does not contain two adjacent elements.
Therefore $C$ is a \1-perfect code by the definition, where $V_{od}=C_1$.
\proofend

\subsection[A Class Of $1'$-Perfect Codes In $Z_{2^k}^n$]{A Class Of \1-Perfect Codes In $Z_{2^k}^n$} \label{ss_b}
Let $n=2^r$ and $I=(i_1,\ldots,i_k)$ be a collection of nonnegative integers such that
$1 i_1+2 i_2+ \ldots +k i_k=r$.
Let $\bar b_1,\ldots,\bar b_n\in Z_{2^k}^{1+i_1+...+i_k}$ be all the elements
of $\{1\}\times (2^{k-1} Z_{2^k})^{i_1}\times (2^{k-2} Z_{2^k})^{i_2}
\times \ldots \times (2^{0} Z_{2^k})^{i_k}$ ordered lexicographically.
And let $B_I$ be the matrix with the columns $\bar b_1,\ldots,\bar b_n$.

\bexam
If $k=3$, $i_1=2$, $i_2=1$, $i_3=0$, then $r=4$, $n=2^r=16$, and
$$
B_I=
\left(%
\begin{array}{c}
1111111111111111 \\
0000000044444444 \\
0000444400004444 \\
0246024602460246 \\
\end{array}
\right).
$$
If $k=3$, $i_1=0$, $i_2=i_3=1$, then $r=5$, $n=2^r=32$, and
$$
B_I=
\left(%
\begin{array}{c}
11111111111111111111111111111111 \\
00000000222222224444444466666666 \\
01234567012345670123456701234567 \\
\end{array}
\right).
$$
If $k=3$, $i_1=i_2=i_3=1$, then $r=6$, $n=2^r=64$, and
$$
B_I=
\left(%
\begin{array}{c}
1111111111111111111111111111111111111111111111111111111111111111 \\
0000000000000000000000000000000044444444444444444444444444444444 \\
0000000022222222444444446666666600000000222222224444444466666666 \\
0123456701234567012345670123456701234567012345670123456701234567 \\
\end{array}
\right).
$$
\eexam

\blemma\label{isPerf}
The linear code
${\cal H}_I\triangleq \{\bar h=(h_1,\ldots,h_n)\in Z_{2^k}^n\,
|\,\sum_{j=1}^n h_j \bar b_j = B_I\bar h^T = \bar 0 \}$ with the check matrix $B_I$
is \1-perfect.
\elemma
\proofr
First we claim that

(*) {\em the distance $2$ between codewords
$\bar h_1$, $\bar h_2$ from ${\cal H}_I$ is impossible}.
Indeed, if $d^\diamond(\bar h_1,\bar h_2)=2$, then the codeword
$\bar h \triangleq \bar h_1-\bar h_2$ have one or two nonzero positions.
The first case contradicts to the fact that $\bar b_i\neq \bar 0$, $i=1, \ldots ,n$.
In the second case we have $\beta \bar b_i+\gamma \bar b_j=\bar 0$ for some different $i$, $j$
and odd $\beta$, $\gamma$. But the first row of $B_I$ implies that $\beta=-\gamma$ and,
since $\beta$ and $\gamma$ are odd, we get $\bar b_i=\bar b_j$.
This again contradicts to the construction of $B_I$. The claim (*) is proved.

First row of 
$B_I$ implies that the words of ${\cal H}_I$ have
even weight.

We need to check that each
odd word is at the distance $1$ from exactly one codeword.
Let $\bar z=(z_1, \ldots ,z_n)\in Z_{2^k}^n$ have an odd weight,
$\bar s\triangleq\sum_{j=1}^n z_i\bar b_i$, and $s\triangleq\sum_{j=1}^n z_i$.
Note that $s$ is odd.
Since $s^{-1}\bar s\in \{1\}\times (2^{k-1} Z_{2^k})^{i_1}\times (2^{k-2} Z_{2^k})^{i_2}
\times \ldots \times (2^{0} Z_{2^k})^{i_k}$, there exists $j'$ such that
$s^{-1}\bar s=\bar b_{j'}$. Let $\bar z'\in Z_{2^k}^n$ be the word with $s$ in the $j'$th position
and zeroes in the others.
It is easy to check that $\bar z-\bar z'$ is a codeword at the distance $1$ from $\bar z$.
As follows from (*), such codeword is unique for each odd $\bar z$.
\proofend

\subsection{Co-$Z_{2^k}$-Linear Extended $1$-Perfect Binary Codes} \label{ss_c}

Now we have all we need to construct a class of co-$Z_{2^k}$-linear
extended $1$-perfect binary codes.
As we will see in \hyperref[S_nonexist]{Section~}\ref{S_nonexist},
the class constructed exhaust all such codes provided the mapping $\Phi$ is fixed.

\btheorem\label{perf}
The co-$Z_{2^k}$-linear code $\tilde H_I\triangleq \Phi({\cal H}_I)$ is
a binary $(nm,2^{nm}/2nm,4)$ code, i.\,e., an extended $1$-perfect code.
\etheorem

\proofr
The statement follows directly from
\hyperref[DcoLin]{Lem\-ma~}\ref{DcoLin},
\hyperref[DcoLin]{Pro\-po\-si\-ti\-on~}\ref{1perf}, and
\hyperref[DcoLin]{Lem\-ma~}\ref{isPerf}.
\proofend

\subsection{Note On The General Case $Z_{2m}$} \label{ss_d}

Indeed, \1-perfect codes can be constructed
over $Z_{2m}$ with $d^\diamond$-distance
for each $m\geq 1$, see the following examples.
Since $m=2^\mu$ is a necessary condition for constructing binary codes
using the way of \hyperref[S_co-Z-lin]{Subsection~}\ref{S_co-Z-lin},
the classification of \1-perfect codes
in the other cases is out of view of this paper.

\bexam\label{ex1perf} The codes with the check matrices

$$B'\triangleq
\left(%
\begin{array}{cccccccc}
  1 & 1 & 1 & 1 & 1 & 1 & 1 & 1  \\
  0 & 0 & 0 & 0 & 12& 12& 12& 12 \\
  0 & 6 & 12& 18& 0 & 6 & 12& 18
\end{array}%
\right),
$$
$$
B''\triangleq
\left(%
\begin{array}{cccccc}
  1 & 1 & 1  & 1  & 1  & 1  \\
  0 & 4 & 8  & 12 & 16 & 20
\end{array}%
\right)
$$
are \1-perfect codes over $Z_{24}$ with $d^\diamond$-distance.
\eexam

\section[$Z_{2^k}$-Linear {H}a\-da\-mard Codes]{%
$Z_{2^k}$-Linear {H}a\-da\-mard Codes}\label{S_Had}

\blemma \label{L_HD}
Let ${\cal H}$ be a linear code in $Z_{2^k}^n$.
Then ${\cal H}$ is an $(n,n2^{k},n2^{k-2})^*$ code
if and only if ${\cal H}^\bot$ is
a \1-perfect code in $(Z_{2^k}^n,d^\diamond)$.
\elemma
\proofr
Assume the code ${\cal H}$ has parameters $(n,n2^{k},n2^{k-2})^*$.
Then $\varphi({\cal H})$ has the parameters of an Hadamard code,
see \hyperref[Cor1]{Corollary~}\ref{Cor1}. By \hyperref[duality]{Theorem~}\ref{duality}
the code $\Phi({\cal H}^\bot)$ is formally dual to the Hadamard code $\varphi({\cal H})$, i.\,e,
$\Phi({\cal H}^\bot)$ is an extended $1$-perfect code. Then the cardinality and
the minimal $d^\diamond$-distance $4$ of the code ${\cal H}^\bot$
follow from \hyperref[DcoLin]{Lemma~}\ref{DcoLin}.

The reverse statement can be proved by reversing the arguments.
\proofend

The next theorem follows immediately from \hyperref[isPerf]{Lemmas~}\ref{isPerf}
\hyperref[L_HD]{and~}\ref{L_HD}.

\btheorem
The code
${\cal D}_I\triangleq {\cal H}_I^\bot$
is a linear $(n,n2^{k},n2^{k-2})^*$ code with the generator matrix $B_I$.
The $Z_{2^k}$-linear code $D_I\triangleq \varphi({\cal D}_I)$ is
a binary $(n2^{k-1},n2^{k},n2^{k-2})$ code, i.\,e., an Ha\-da\-mard code.
\etheorem

\bnote
The code ${\cal D}_{(r,0,\ldots,0)}$ is the
{\em first order Reed-Muller code over $Z_{2^k}$} \cite{GuptaOth2005}.
\enote

\bnote Indeed, if an Hadamard code
$A$ of length $m$ exists (see \hyperref[S_Z-lin0]{Subsection~}\ref{S_Z-lin0}),
then $Z_{2m}$-linear Hadamard code can be constructed for length $nm=2^r m$,
even if $m$ is not a power of two.
For example, the linear $Z_{24}$-code with the generator matrix $B'$
from  \hyperref[ex1perf]{Example~}\ref{ex1perf}
has parameters $(16,192,48)^*$ and the corresponding binary $Z_{24}$-linear
code has parameters $(96,192,48)$, i.\,e., the parameters of an Hadamard code of length $96$.

However, the code with the generator matrix $B''$ (\hyperref[ex1perf]{Example~}\ref{ex1perf})
has parameters $(6,144,30)^*$, corresponding a binary $(72,144,30)$ code.
The code distance of this code is smaller than the code distance of an Hadamard code
with the same length and cardinality.
\enote

\section[Nonexistence Of Unknown Co-$Z_{2^k}$-Linear Extended Perfect Codes And
$Z_{2^k}$-Linear {H}a\-da\-mard Codes]{%
Nonexistence Of Unknown Co-$Z_{2^k}$-Linear Extended Perfect Codes And
$Z_{2^k}$-Linear {H}a\-da\-mard Codes}\label{S_nonexist}
{L}{et} $n$ be a power of $2$. In this section we will show
(\hyperref[Th_HH]{Theorem~}\ref{Th_HH}) that
each linear $(n, 2^{kn}/n2^k, 4)^\diamond$ code in $Z_{2^k}^n$ is equivalent to a code
from the class constructed in \hyperref[S_perf]{Section~}\ref{S_perf}. Similarly,
each linear
$(n,n2^{k},n2^{k-2})^*$
 code in $Z_{2^k}^n$ is equivalent to a code
from the class constructed in \hyperref[S_Had]{Section~}\ref{S_Had}
(\hyperref[Th_DD]{Theorem~}\ref{Th_DD}). The key moment of the proof
is \hyperref[Th_HH]{Lemma~}\ref{L_1111}.
The partial $Z_4$ case of this lemma
was proved in \cite{Kro:2000:Z4_Perf} and in \cite{BorPheRif:2003},
but the proof given here is not similar.

We say that two linear codes ${\cal C}_1,{\cal C}_2\subseteq Z_{2^k}^n$
are \emph{equivalent} if 
${\cal C}_2=\bar z\circ\pi {\cal C}_1$
where $\pi$ is a coordinate permutation,
$\bar z \in (Z_{2^k}^*)^n\triangleq \{1,3,\ldots,2^k-1\}^n$,
and $\circ$ is the coordinate-wise product defined as
$(z_1,\ldots,z_n)\circ(x_1,\ldots,x_n)\triangleq (z_1 x_1,\ldots,z_n x_n)$.
Note that both operations $\pi$ and $\bar z\,\circ$ are group automorphisms of the
additive group of $Z_{2^k}^n$ and isometries of the metric spaces
$(Z_{2^k}^n,d^*)$ and $(Z_{2^k}^n,d^\diamond)$.
The following proposition shows that there are no other linear isometries
of $(Z_{2^k}^n,d^*)$ or $(Z_{2^k}^n,d^\diamond)$, and thus our definition
of the equivalence is natural.
\bpro
Assume that $\Gamma$ is a linear transformation of $Z_{2^k}^n$ and
isometry of $(Z_{2^k}^n,d)$ where $d=d^*$ or $d=d^\diamond$.
Then $\Gamma(\bar x)\equiv\bar z\circ\pi(\bar x)$
where $\pi$ is a coordinate permutation,
and $\bar z \in (Z_{2^k}^*)^n$.
\epro
\proofr
Denote by $\bar e_i$ the word with $1$ in $i$th
position and zeroes in the other positions.
It is enough to check that for each $i$ we have
$\Gamma(\bar e_i)\equiv z_j \bar e_{j}$ with some
$j=\pi(i)$ and $z_{j}\in Z_{2^k}^*$.
Indeed, from the isometric properties of $\Gamma$
we derive that $\Gamma(\bar e_i)$ has only one non-zero
coordinate.
On the other hand, this coordinate belongs to $Z_{2^k}^*$,
because $\Gamma$ is a bijection.
\proofend

The following two theorems are the main results of this section.
Remind that the codes ${\cal H}_I$ and ${\cal D}_I$ are defined in
\hyperref[S_perf]{Sections~}\ref{S_perf} \hyperref[S_Had]{and~}\ref{S_Had}.

\btheorem \label{Th_HH}
Let ${\cal H}\subset Z_{2^k}^n$ be a linear
\1-perfect
code.
Then ${\cal H}$ is equivalent to ${\cal H}_I$ where
$I=(i_1,\ldots,i_k)$ is a collection of nonnegative integers such that
$1 i_1+2 i_2+ \ldots +k i_k=\log_2 n$.
\etheorem

\btheorem \label{Th_DD}
Let ${\cal D}\subset Z_{2^k}^n$ be a linear $(n,n2^{k},n2^{k-2})^*$ code.
Then ${\cal D}$ is equivalent to ${\cal D}_I$, where
$I=(i_1,\ldots,i_k)$ is a collection of nonnegative integers such that
$1 i_1+2 i_2+ \ldots +k i_k=\log_2 n$.
\etheorem

By \hyperref[L_HD]{Lemma~}\ref{L_HD}
it is enough to prove only one of \hyperref[Th_HH]{Theorems~}\ref{Th_HH}, \ref{Th_DD}.
We need the following auxiliary result.

\blemma\label{L_1111}
If $\cal D$ is a linear $(n,n2^{k},n2^{k-2})^*$ code in $Z_{2^k}^n$, then
$\cal D$ contains an element from $(Z_{2^k}^*)^n$.
\elemma
\proofr
We prove the lemma by induction.
If $n=1$,
then ${\cal D}=Z_{2^k}$.
Assume $n>1$.

(*) \emph{We claim that $\cal D$ contains an element from $\{0,2^{k-1}\}^n$ of weight $n2^{k-2}$.}
Let $\bar x'$ and $\bar x''$ be two linear independent elements in $\cal D$
of order $2m'$ and $2m''$ respectively.
This means that $m' \bar x'$ and $m'' \bar x''$
are different nonzero elements from ${\cal D} \cap \{0,2^{k-1}\}^n$.
Since $\varphi$ is an isometric imbedding (see \hyperref[P_iso]{Proposition~}\ref{P_iso})
and $\varphi({\cal D})$ is an $(n2^{k-1},n2^{k},n2^{k-2})$ Hadamard code,
the only possible values of
$wt^*$-weight of elements in $\cal D$ are $0$, $n2^{k-2}$ and $n2^{k-1}$.
So, at least one of $m' \bar x'$ and $m'' \bar x''$ has weight $n2^{k-2}$.
The claim (*) is proved.

Without loss of generality
assume that
\(\bar a \triangleq (2^{k-1},\ldots,2^{k-1},0,\ldots,0)
\in {\cal D}.\)

Let us consider two codes obtained from $\cal D$ by puncturing $n/2$ coordinates:
\begin{IEEEeqnarray*}{rCl}
{\cal D}_1 & {}\triangleq{} & \{\bar z'\in Z_{2^k}^{n/2} \, | \  %
\exists\, \bar z''\in Z_{2^k}^{n/2} :\, (\bar z',\bar z'')\in {\cal D}\}, \\
{\cal D}_2 & {}\triangleq{} & \{\bar z''\in Z_{2^k}^{n/2} \, | \  %
\exists\, \bar z'\in Z_{2^k}^{n/2} :\, (\bar z',\bar z'')\in {\cal D}\}.
\end{IEEEeqnarray*}

(**) \emph{We claim that ${\cal D}_1$ and ${\cal D}_2$ are linear $(n/2,n2^{k-1},n2^{k-3})^*$ codes.}
The linearity of the codes is obvious.
The code distance follows from the fact that the distance between $\bar a$
and any element $\bar z\in {\cal D}$ is $0$, $n2^{k-2}$ or $n2^{k-1}$.
It is true that $|{\cal D}_1| \geq |{\cal D}|/2$ because the code distance of ${\cal D}$
permits of only one nonzero codeword with zeroes in the first $n/2$ coordinates.
On the other hand, such a codeword exists:
$(0,\ldots,0,2^{k-1},\ldots,2^{k-1})=(2^{k-1},\ldots,2^{k-1})+\bar a\in {\cal D}$.
So, $|{\cal D}_1| = |{\cal D}|/2 = n2^{k-1}$.
Similarly, we get $|{\cal D}_2| = n2^{k-1}$. The claim (**) is proved.

By the assumption of induction ${\cal D}_1$ contains an element $\bar u'$
from $(Z_{2^k}^*)^{n/2}$.
This means that there exist $\bar u''$ from $Z_{2^k}^{n/2}$ such that
$(\bar u',\bar u'')\in {\cal D}$.
Since $wt^*(2^{k-1}(\bar u',\bar u''))\in \{n2^{k-2},n2^{k-1}\}$ and
$2^{k-1}\bar u'=(2^{k-1},...,2^{k-1})$,
we have $wt^*(2^{k-1}\bar u'')=0$ or $wt^*(2^{k-1}\bar u'')=n2^{k-2}$.
So,
\setlength{\arraycolsep}{0mm}
\begin{eqnarray}
 \label{u''0} 2^{k-1}\bar u'' & {}={} & (0,\ldots,0), \\
 \label{u''1} \mbox{or}\ \ 2^{k-1}\bar u'' & {}={} & (2^{k-1},\ldots,2^{k-1}).
\end{eqnarray}
Similarly, considering the code ${\cal D}_2$, we can
find a codeword $(\bar v',\bar v'')\in {\cal D}$ such that
$2^{k-1}\bar v''=(2^{k-1},...,2^{k-1})$ and $v'$ satisfies\\[-20pt]
\setlength{\arraycolsep}{0mm}
\begin{eqnarray}
 \label{v'0} 2^{k-1}\bar v' & {}={} & (0,\ldots,0), \\
 \label{v'1} \mbox{or}\ \ 2^{k-1}\bar v' & {}={} & (2^{k-1},\ldots,2^{k-1}).
\end{eqnarray}
If (\ref{u''1}) is true, then $(u',u'')\in (Z_{2^k}^*)^{n}\cap {\cal D}$.
If (\ref{v'1}) is true, then $(v',v'')\in (Z_{2^k}^*)^{n}\cap {\cal D}$.
If (\ref{u''0}) is true and (\ref{v'0}) is true, then
$(u',u'')+(v',v'')\in (Z_{2^k}^*)^{n}\cap {\cal D}$.
\hyperref[L_1111]{Lemma~}\ref{L_1111} is proved.
\proofend

\proofrem{Proof of \hyperref[Th_DD]{Theorem~}\ref{Th_DD}}
By \hyperref[L_1111]{Lemma~}\ref{L_1111} the code
$\cal D$ contains an element $\bar c=(c_1,\ldots,c_n)$ from $(Z_{2^k}^*)^n$.
Then the code ${\cal D}'\triangleq (c_1^{-1},\ldots,c_n^{-1})\circ {\cal D}$
is equivalent to ${\cal D}$ and contains $\bar 1=(1,\ldots,1)$.
Let $\{\bar 1, q_1,\ldots, q_s\}$ be a basis of ${\cal D}'$ and
$i_j$ be the number of elements of order $2^{j}$, $j=1,\ldots,k$.
Assume $\bar 1, q_1,\ldots, q_s$ are the rows of the matrix $Q$
of size $(1+i_1+\ldots+i_k)\times n$, a generator matrix of the code ${\cal D}'$.
Then the columns of $Q$ are elements
of $\{1\}\times (2^{k-1} Z_{2^k})^{i_1}\times (2^{k-2} Z_{2^k})^{i_2}
\times \ldots \times (2^{0} Z_{2^k})^{i_k}$.
Since by \hyperref[L_HD]{Lemma~}\ref{L_HD} the code $d^\diamond$-distance of ${\cal D}'^\bot$ is more than $2$,
all the columns are pairwise different.
Therefore $Q$ coincides with $B_I$, $I=(i_1,\ldots,i_k)$,
up to permutation of columns.
\proofremend

So, we conclude, provided the mappings $\varphi$ and $\Phi$ are fixed,
all up to equivalence
co-$Z_{2^k}$-linear extended $1$-perfect codes and
$Z_{2^k}$-linear Hadamard codes are described
in \hyperref[S_perf]{Sections~}\ref{S_perf} \hyperref[S_Had]{and~}\ref{S_Had}.

\ifx\href\undefined \newcommand\href[2]{#2} \fi\ifx\url\undefined
  \newcommand\url[1]{\href{#1}{#1}} \fi

\end{document}


\renewcommand\today{}

\title{О $Z_{2^k}$-дуальных двоичных кодах%
\thanks{ Это авторский перевод статьи в
IEEE Transactions on Information Theory 53(4) 2007, 1532--1537, Digital
Object Identifier \href{http://dx.doi.org/10.1109/TIT.2007.892787}{10.1109/TIT.2007.892787}, \copyright
2007 IEEE.}
\thanks{ Результаты работы частично докладывались на 4й Международной конференции по оптимальным кодам и смежным вопросам OC\,2005 (Пампорово, Болгария, июнь 2005).}
}
\renewcommand\footnotemark{ }
\author{Д.~С.~Кротов%
\thanks{ Адрес автора: Институт математики им. С.\,Л.\,Соболева СО РАН, проспект Академика Коптюга 4, Новосибирск 630090, Россия
(e-mail: {\tt krotov@math.nsc.ru}).}
}

\markboth
{Krotov: On $Z_{2^k}$-Dual Binary Codes}
{Krotov: On $Z_{2^k}$-Dual Binary Codes}

\maketitle

\begin{abstract}
Представлено новое обобщение отображения Грея $\Phi: Z_{2^k}^n \to Z_{2}^{2^{k-1}n}$,
связанное с известным обобщенным отображением Грея $\varphi$ следующим образом:
если взять два дуальных линейных $Z_{2^k}$-кода и построить из них двоичные коды, используя обобщения $\varphi$ и $\Phi$ отображения Грея, то весовые энумераторы полученных двоичных кодов будут связаны тождеством Мак-Вильямс.
Описаны классы $Z_{2^k}$"=линейных кодов Адамара и ко-$Z_{2^k}$"=линейных расширенных 1"=совершенных кодов, где ко-$Z_{2^k}$"=линейность означает, что код может быть получен из линейного $Z_{2^k}$-кода при помощи нового обобщенного отображения Грея.

Ключевые слова: отображение Грея, код Адамара, тождество Мак-Вильямс, совершенный код, $Z_{2^k}$"=линейность
\end{abstract}

\section{Введение}\label{cl_ru:S_Intro}
Как обнаружено в \cite{cl_ru:Nechaev82}, \cite{cl_ru:Nechaev:Kerdock_cyclic_form},
некоторые нелинейные двоичные коды представимы как линейные коды над кольцом $Z_4$.
Вариант такого представления, найденный в \cite{cl_ru:HammonsOth:Z4_linearity},
использует \emph{отображение Грея} $\phi: 0{\to} 00,\ 1{\to} 01,\ 2{\to} 11,\ 3{\to} 10$ для построения двоичных так называемых $Z_4$"=линейных кодов из кодов в четырехбуквенном алфавите.
Ключевое свойство отображения $\phi$ с этой точки зрения "--- что оно является изометрией между $Z_4$ с метрикой Ли и $Z_2^2$ с метрикой Хемминга.
В работе \cite{cl_ru:Carlet:Z2k-linear} (и, в более общей форме, в \cite{cl_ru:HonNech})
отображение Грея обобщено для построения $Z_{2^k}$"=линейных кодов.
Обобщенное отображение Грея (скажем, $\varphi$; см. \hyperref[cl_ru:S_Z-lin]{Подраздел~}\ref{cl_ru:S_Z-lin} для напоминания основных фактов об обобщенном отображении Грея)
является изометричным вложением $Z_{2^k}$ с метрикой, определяемой однородной весовой функцией \cite{cl_ru:ConstHeise}, в $Z_2^{2^{k-1}}$ с метрикой Хемминга.

В настоящей статье мы представляем другое обобщение $\Phi$ отображения Грея (\hyperref[cl_ru:S_co-Z-lin]{Подраздел~}\ref{cl_ru:S_co-Z-lin}).
Это отображение оказалось дуальным предыдущему в следующем смысле.
Если $\cal C$ и ${\cal C}^\bot$ "--- дуальные линейные $Z_{2^k}$-коды,
то двоичный $Z_{2^k}$"=линейный код $\varphi({\cal C})$ и {\it ко-$Z_{2^k}$"=линейный} код $\Phi({\cal C}^\bot)$ формально дуальны.
Формальная дуальность означает, что весовые энумераторы этих двух кодов удовлетворяют тождеству Мак-Вильямс (\hyperref[cl_ru:S_dual]{Раз\-дел~}\ref{cl_ru:S_dual}) (заметим, что эти коды в общем случае могут быть нелинейными, и тогда они не могут быть дуальными в обычном смысле, как подпространства двоичного векторного пространства).
Таким образом, мы решаем проблему дуальности для $Z_{2^k}$"=линейных двоичных кодов: связать весовые энумераторы двоичных образов дуальных линейных кодов над $Z_{2^k}$. Эта проблема не могла быть решена посредством только стандартного обобщенного отображения Грея $\varphi$, поскольку, как замечено в \cite{cl_ru:Carlet:Z2k-linear},
весовые энумераторы кодов $\varphi({\cal C})$ и $\varphi({\cal C}^\bot)$ в общем случае не являются связанными, в отличие от случая $Z_4$"=линейных кодов.

В \cite{cl_ru:HonNech} показано, что двоичные (и не только двоичные) коды могут быть представлены как групповые коды над различными группами.
Такие представления используют особую \emph{масштабную изометрию}
($\varphi$ "--- частный случай такой изометрии),
действующую изометрично из некоторого модуля со специально определенной метрикой в двоичное (или, в общем случае, не только двоичное) пространство Хемминга.
С таким подходом каждый элемент модуля соответствует одному кодовому слову, $Z_{2^k}$"=линейные коды являются частным случаем подобного представления.

В нашем подходе каждый элемент модуля (слово) соответствует некоторому множеству кодовых слов, являющемуся декартовым произведением множеств, соответствующих (посредством специального отображения) символам слова модуля.
В конструктивной части этот подход может быть представлен обобщенной каскадной конструкцией \cite{cl_ru:Zin1976:GCC}, но кодовое расстояние полученного кода отличается от конструктивного расстояния обобщенного каскадного кода.
Вычисляя расстояние, мы используем определенные изометрические свойства отображения, которое может быть рассмотрено как некоторый аналог масштабной изометрии.
Отображение $\Phi$, рассматриваемое в настоящей работе, позволяет строить только коды с расстоянием не больше $4$, но общий подход может иметь больший потенциал.

Новый способ обобщения отображения Грея (\hyperref[cl_ru:S_co-Z-lin]{Подраздел~}\ref{cl_ru:S_co-Z-lin})
и его дуальность (\hyperref[cl_ru:S_dual]{Раздел~}\ref{cl_ru:S_dual}) <<старому способу>> "--- основные результаты первой части работы.
Вторая часть (\hyperref[cl_ru:S_perf]{Разделы~}\ref{cl_ru:S_perf} \hyperref[cl_ru:S_Had]{--} \ref{cl_ru:S_nonexist})
посвящена построению ко-$Z_{2^k}$"=линейных и $Z_{2^k}$"=линейных кодов с параметрами расширенных $1$"=совершенных кодов и кодов Адамара.
В серии статей
(напр.
\cite{cl_ru:RifPuj:1997},
\cite{cl_ru:BorRif:1999},
\cite{cl_ru:Kro:2000:Z4_Perf},
\cite{cl_ru:Kro:2001:Z4_Had_Perf},
\cite{cl_ru:PheRif:2002},
\cite{cl_ru:BorPheRif:2003},
\cite{cl_ru:PheRifVil:2004ACCT:additive_hadamard},
\cite{cl_ru:PheRifVil:Hadamard-rank-kernel},
\cite{cl_ru:GuptaOth2005})
рассматриваются $1$"=совершенные (или расширенные $1$"=совершенные) коды и коды Адамара, являющиеся линейными в некотором нестандартном смысле (отметим, что для каждого из рассматриваемых значений кодовых параметров существует только один, с точностью до эквивалентности, линейный код).
Пользуясь двумя обобщенными отображениями Грея $\varphi$ и $\Phi$, мы строим широкий класс таких кодов.
Некоторые естественные вопросы о построенных кодах остаются открытыми для дальнейшего исследования: какие из этих кодов эквивалентны друг другу или другим известным кодам; какие из них являются пропелинейными (см. определение в \cite{cl_ru:RifPuj:1997}) или по крайней мере транзитивными
(см. определение и некоторые конструкции транзитивных $1$"=совершенных кодов в \cite{cl_ru:Sol:2005transitive});
установить границы размерности ядра и ранга этих кодов (размерность ядра и ранг "--- хорошие меры линейности нелинейных кодов в классическом двоичном случае;
для $Z_4$"=линейных расширенных $1$"=совершенных кодов и кодов Адамара эти параметры вычислены в работах \cite{cl_ru:Kro:2000:Z4_Perf}, \cite{cl_ru:Kro:2001:Z4_Had_Perf}, \cite{cl_ru:BorPheRif:2003}, \cite{cl_ru:PheRifVil:2004ACCT:additive_hadamard}, \cite{cl_ru:PheRifVil:Hadamard-rank-kernel});
и т.\,д.

В \hyperref[cl_ru:S_perf]{Разделе~}\ref{cl_ru:S_perf} мы строим
ко-$Z_{2^k}$"=линейные расширенные $1$"=совершенные $(n,2^n/2n,4)$-коды.
Как отмечено выше, кодовое расстояние ко-$Z_{2^k}$"=линейного кода
не может превышать $4$ при $k>2$ (см. \hyperref[cl_ru:DcoLin]{Лемму~}\ref{cl_ru:DcoLin}).
Таким образом, расширенный $1$"=совершенные коды "--- наилучшие примеры.

Как следствие, дуальное расстояние $Z_{2^k}$"=линейных кодов ($k>2$) также не превышает $4$.
Наилучшие примеры "--- коды с большим кодовым расстоянием.
В \hyperref[cl_ru:S_Had]{Разделе~}\ref{cl_ru:S_Had} мы строим $Z_{2^k}$"=линейные коды с параметрами $(n,2n,n/2)$, т.\,е. коды Адамара.
Они $Z_{2^k}$-дуальны расширенным $1$"=совершенным кодам из \hyperref[cl_ru:S_perf]{Раздела~}\ref{cl_ru:S_perf}, т.\,е., $\varphi$- и, соответственно, $\Phi$- прообразы построенных расширенных $1$ совершенных кодов и кодов Адамара являются дуальными $Z_{2^k}$-кодами.

Описанные конструкции обобщают конструкции $Z_4$"=линейных расширенных $1$"=совершенных кодов и кодов Адамара \cite{cl_ru:Kro:2000:Z4_Perf}, \cite{cl_ru:Kro:2001:Z4_Had_Perf} (см. также \cite{cl_ru:BorPheRif:2003}).

Естественный вопрос "--- является ли приведенная конструкция ко-$Z_{2^k}$"=линейных расширенных $1$"=совершенных кодов и $Z_{2^k}$"=линейных кодов Адамара полной или нет.
В \hyperref[cl_ru:S_nonexist]{Разделе~}\ref{cl_ru:S_nonexist} мы показываем, что ответ положительный: все ко-$Z_{2^k}$"=линейные расширенные $1$"=совершенные коды и все $Z_{2^k}$"=линейные коды Адамара эквивалентны кодам, построенным в \hyperref[cl_ru:S_perf]{Разделах~}\ref{cl_ru:S_perf} \hyperref[cl_ru:S_Had]{и~}\ref{cl_ru:S_Had}.

\section[Определения и основные факты]{Определения и основные факты. Два обобщения отображения Грея}\label{cl_ru:S_def}
Мы будем использовать следующие общие обозначения.
\emph{Двоичный $(n,M,d)$-код} "--- подмножество мощности $M$ множества $Z_2^n = {Z_2\times\ldots\times Z_2}$ с расстоянием не меньше $d$ между каждыми двумя различными элементами (\emph{кодовыми словами});
$n$ называется \emph{длиной} кода, $d$ "--- \emph{кодовым расстоянием}.
Двоичный $(n,2n,n/2)$-код называется \emph{кодом Адамара}; такие коды существуют если и только если существуют матрицы Адамара размера $n\times n$, т.\,е. как минимум $n\equiv 0\mbox{ mod }4$ при $n>2$, см. напр. \cite{cl_ru:MWS}.
Линейный код Адамара известен также как код Рида-Маллера первого порядка, в этом случае $n$ "--- степень двойки.
Двоичные $(n,2^n/2n,4)$-коды называются \emph{расширенными $1$"=совершенными кодами}; они существуют если и только если $n$ является степенью двойки.
Такой код всегда является расширением некоторого $(n{-}1,2^n/2n,3)$-кода, который называется $1$"=совершенным, таким образом, изучение расширенных $1$"=совершенных двоичных кодов суть альтернативный подход к изучению $1$"=совершенных двоичных кодов. Для каждой допустимой длины существует только один, с точностью до перестановки координат, линейный (расширенный) $1$"=совершенный код и, при $n>8$, много нелинейных, классификация которых в настоящее время является открытой проблемой.
В следующих подразделах мы определим основные понятия, которые будут использоваться в статье: две весовые функции в $Z_{2^k}$, две соответствующие метрики, два обобщенных отображения Грея $\varphi$ и $\Phi$, понятия $Z_{2^k}$"=линейных и ко-$Z_{2^k}$"=линейных кодов, сформулируем простые но фундаментальные утверждения  об изометрических свойствах отображений $\varphi$ и $\Phi$.
\subsection[$Z_{2m}$-Линейные коды]{$Z_{2m}$-Линейные коды}\label{cl_ru:S_Z-lin}

В этом подразделе мы напомним основные определения и факты об обобщенном отображении Грея и $Z_{2^k}$"=линейных кодах \cite{cl_ru:Carlet:Z2k-linear},  \cite{cl_ru:HonNech}.

Пусть $m \geq 2$ "--- такое целое число, при котором существует матрица Адамара размера $m\times m$.
Пусть $A\subset Z_2^{m}$ есть $(m,2m,m/2)$-код Адамара и $A=\{a_0,a_1, \ldots ,a_{2m-1}\}$, где $a_0$ есть слово из всех нулей и $a_i+a_{i+m}$ есть слово из всех единиц для каждого $i$ от $0$ до $m-1$.
Определим {\em обобщенное отображение Грея} $\varphi : Z_{2m}^n \to Z_{2}^{mn}$ следующим правилом:
$$
\varphi(x_1,\ldots,x_n)\triangleq (a_{x_1},\ldots,a_{x_n}).
$$
Весовую функцию $wt^*:Z_{2m}\to R^+$ определим равенством
$$
wt^*(x) \triangleq
\cases{
0& \mbox{если }$x=0$, \cr
m& \mbox{если }$x=m$, \cr
m/2& \mbox{если }$x\neq 0,m$.}
$$
В случае, когда $m$ "--- степень двойки, функция $wt^*$ является \emph{однородным весом}, представленным в \cite{cl_ru:ConstHeise} для более общего класса колец.
Соответствующее расстояние $d^*$ в $Z_{2m}^n$ определяется стандартным способом:
$$ d^*((x_1,\ldots,x_n),(y_1,\ldots,y_n))
\triangleq
\sum_{i=1}^n wt^*(y_i-x_i).
$$

\begin{pro}[\cite{cl_ru:HonNech}] \label{cl_ru:P_iso}
Отображение $\varphi$ является изометричным вложением $(Z_{2m}^n,d^*)$ в
$(Z_2^{mn},d)$, т.\,е.
$d^*(\bar x,\bar y)=d(\varphi(\bar x),\varphi(\bar y))$, где $d(\cdot,\cdot)$ "--- расстояние Хемминга.
\end{pro}

Если ${\cal C}\subseteq Z_{2m}^n$, то
$$
\varphi({\cal C})\triangleq \{ \varphi(\bar x) \,|\, \bar x\in {\cal C} \}.
$$

Будем говорить, что ${\cal C}$ есть $(n,M,d)^*$-код, если 
${\cal C}\subseteq Z_{2m}^n$, $|{\cal C}|=M$ и $d^*$-расстояние между любыми двумя различными элементами ${\cal C}$ не меньше $d$.

\begin{corol}\label{cl_ru:Cor1}
Если ${\cal C}$ является $(n,M,d)^*$-кодом, то $\varphi({\cal C})$ "--- двоичный $(mn,M,d)$-код.
\end{corol}

Двоичный код называется {\em $Z_{2m}$"=линейным}, если его координаты могут быть упорядочены таким образом, что он является образом некоторого линейного $Z_{2m}$-кода под действием отображения $\varphi$.
Как отмечено в \cite{cl_ru:Carlet:Z2k-linear}, длина $Z_{2m}$"=линейного кода должна делиться на $m$, а вес каждого кодового слова должен делиться на $m/2$.

\begin{note}
Отображение $\varphi$ и понятие $Z_{2m}$"=линейности, приведенное выше, зависит от выбора кода Адамара $A$ и единственное ограничение на $m$ "--- существование матрицы Адамара размера $m\times m$.
Раньше {\rm \cite{cl_ru:Carlet:Z2k-linear}}
(и в двоичном подслучае в \cite[Раздел D]{cl_ru:HonNech}) рассматривался только случай
$m=2^{k-1}$ и в качестве $A$ брался код $R(1,k-1)$, код Рида-Маллера первого порядка.
\end{note}

\subsection[Ко-$Z_{2^k}$"=линейные коды]{Ко-$Z_{2^k}$"=линейные коды}
\label{cl_ru:S_co-Z-lin}
В этом подразделе мы предложим другой подход к построению двоичных кодов из линейных кодов над $Z_{2^k}$.

Сначала мы представим новый способ обобщить отображение Грея, определив отображение
$\Phi : Z_{2m}^n \to 2^{Z_{2}^{mn}}$,
где $2^{Z_{2}^{mn}}$ обозначает множество всех подмножеств $Z_{2}^{mn}$.
Затем мы определим весовую функцию $wt^\diamond:Z_{2m}\to R^+$
и соответствующее расстояние
$d^\diamond$ и установим (\hyperref[cl_ru:DcoLin]{Лемма~}\ref{cl_ru:DcoLin})
связь между $d^\diamond$-расстоянием кода ${\cal C}\subset Z_{2m}^n$
и расстоянием Хемминга кода $\Phi({\cal C})$.
Наконец, мы определим понятие ко-$Z_{2^k}$"=линейного кода.

Положим $m=2^{k-1}$. Пусть $\{ H_0,\ldots,H_{2m-1} \}$ "--- разбиение $Z_2^{m}$ на
расширенные $1$"=совершенные $(m,2^m/2m,4)$-коды
(например, в качестве $H_0$ мы можем взять расширенный код Хемминга а в качестве остальных частей "--- смежные классы по нему).
Более того, будем полагать, что $H_0$ содержит слово из всех нулей $\bar 0$ и
четность весов кодовых слов из $H_j$ совпадает с четностью $j$.

Определим отображение $\Phi : Z_{2m}^n \to 2^{Z_{2}^{mn}}$ по правилу
$$
\Phi(x_1,\ldots,x_n)\triangleq H_{x_1}\times\ldots\times H_{x_n}.
$$
Если ${\cal C}\subseteq Z_{2m}^n$, то\\[-10pt]
$$
\Phi({\cal C})\triangleq \bigcup_{\bar x\in {\cal C}} \Phi(\bar x).
$$

\begin{examp}
Пусть $k=3$, $H_0=\{0000,1111\}$, \ldots,
$H_2=\{1100,0011\}$, \ldots,
 $H_7=\{0001,1110\}$. Тогда
$$\begin{array}{rcrcl}
&&\Phi(2\,0\,7)&=&H_2\times H_0\times H_7\\
&&&=&\{\, 1100\,0000\,0001,\ 1100\,0000\,1110,\ 1100\,1111\,0001,\ 1100\,1111\,1110,\\
&&&&\phantom{\{\,} 0011\,0000\,0001,\ 0011\,0000\,1110,\ 0011\,1111\,0001,\ 0011\,1111\,1110\,\}.
\end{array}$$
\end{examp}
\hspace{1.0ex}

Определим весовую функцию $wt^\diamond:Z_{2m}\to R^+$ следующим образом:
$$
wt^\diamond(x)\triangleq
\cases{
0&\mbox{если $x=0$},\cr
1&\mbox{если $x$ нечетно},\cr
2&\mbox{если $x\neq 0$ четно}.}
$$
Соответствующее расстояние $d^\diamond$ в $Z_{2m}^n$ определяется стандартным образом:
$$ d^\diamond((x_1,\ldots,x_n),(y_1,\ldots,y_n))
\triangleq
\sum_{i=1}^n wt^\diamond(y_i-x_i).
$$

Будем говорить, что ${\cal C}$ есть $(n,M,d)^\diamond$-код, если
${\cal C}\subseteq Z_{2m}^n$,
$|{\cal C}|=M$ и $d^\diamond$-расстояние между
любыми двумя различными элементами из ${\cal C}$ не меньше $d$.

\begin{lemma}\label{cl_ru:DcoLin}
Пусть $m\geq 4$.
Если ${\cal C}\subseteq Z_{2m}^n$ есть $(n,M,d)^\diamond$-код,
то $\Phi({\cal C})$ "--- двоичный $(mn,M (\frac{2^m}{2m})^n,{\min}(4,d))$-код.
\end{lemma}
Доказательство состоит в непосредственной проверке, и мы опустим его.
%
%
%
%
%
%

Назовем двоичный код {\it ко-$Z_{2m}$"=линейным}, если 
его координаты могут быть упорядочены таким образом, чтобы он являлся образом некоторого линейного $Z_{2m}$-кода при отображении $\Phi$.

\begin{note}
Ко-$Z_{2m}$"=линейные коды могут быть рассмотрены как частный случай обобщенных каскадных кодов \cite{cl_ru:Zin1976:GCC}.
Для этого особого случая кодовое расстояние согласно \hyperref[cl_ru:DcoLin]{Леммы~}\ref{cl_ru:DcoLin} может быть больше расстояния, которое гарантирует обобщенная каскадная конструкция.
\end{note}

\section[$Z_{2^k}$-Дуальность двоичных кодов]{%
$Z_{2^k}$-Дуальность двоичных кодов}
\label{cl_ru:S_dual}
{В} этом разделе мы покажем (\hyperref[cl_ru:duality]{Теорема~}\ref{cl_ru:duality}),
что если два кода получены из дуальных $Z_{2^k}$-кодов посредством обобщений $\varphi$ (один код) и $\Phi$ (другой) отображения Грея,
то эти коды являются формально дуальными, т.\,е. их весовые энумераторы связаны тождеством Мак-Вильямс.

Пусть $C$ есть $Z_{2^k}$"=линейный код, $\varphi$-образ
линейного $Z_{2^k}$-кода $\cal C$ длины $n$.
Пусть ${\cal C}^\bot$ "--- линейный $Z_{2^k}$-код, дуальный коду ${\cal C}$.
Пусть $SW_{{\cal C}^\bot}(X,Z,T)$ "--- многочлен, полученный из полного весового энумератора $W_{{\cal C}^\bot}(X_0,X_1,\ldots,X_{2^k-1})$ кода ${\cal C}^\bot$ отождествлением в $Z$ (соответственно в $T$)
всех $X_j$ с нечетным $j$ (соответственно, с четным $j$, отличным от $0$).

\begin{lemma}[\cite{cl_ru:Carlet:Z2k-linear}]\label{cl_ru:CC}
Выполнено тождество
$$\begin{array}{rrl}
W_C(X,Y)&
=&\displaystyle\frac{1}{|{\cal C}^\bot|}SW_{{\cal C}^\bot} \big(\
X^{2^{k-1}}\hspace{-0.8ex}+Y^{2^{k-1}}\hspace{-0.8ex}+(2^k-2)(X Y)^{2^{k-2}}\hspace{-0.4ex},\  \\
&&\phantom{\frac{1}{|{\cal C}^\bot|}SW_{{\cal C}^\bot} \big(\ }
X^{2^{k-1}}\hspace{-0.8ex}-Y^{2^{k-1}},\   \\
&&\phantom{\frac{1}{|{\cal C}^\bot|}SW_{{\cal C}^\bot} \big(\ }
X^{2^{k-1}}\hspace{-0.8ex}+Y^{2^{k-1}}\hspace{-0.8ex}-2(X Y)^{2^{k-2}}\ \big). \\
\end{array}$$
\end{lemma}

Следующая лемма "--- известный факт о весовом распределении расширенных $1$"=совершенных кодов.
\begin{lemma}\label{cl_ru:W_H}
Пусть $H$ "--- расширенный $1$"=совершенный $(m,2^m/2m,4)$-код.
Тогда \\
{\rm a)}
\(
\displaystyle\frac{1}{|H|}
W_H(X+Y,X-Y)=X^{m}+Y^{m}+(2m-2)(X Y)^{m/2}
\) \ если $\bar 0\in H$,
\\[4pt]
{\rm b)}
\(
\displaystyle\frac{1}{|H|}
W_H(X+Y,X-Y)=X^{m}-Y^{m}
\) \ если $H$ нечетновесовой,
\\[4pt]
{\rm c)}
\(
\displaystyle\frac{1}{|H|}
W_H(X+Y,X-Y)=X^{m}+Y^{m}-2(X Y)^{m/2}
\) \ если $H$ четновесовой и $\bar 0\not\in H$.%
\end{lemma}

\begin{pro}\label{cl_ru:W-SW}
Пусть $\cal C$ есть $Z_{2^k}$-код и
$\tilde C=\Phi({\cal C})$. Тогда
$$
W_{\tilde C}(X,Y)=SW_{\cal C}(W_{H_0}(X,Y), W_{H_1}(X,Y), W_{H_2}(X,Y)),
$$
где $\{ H_0,\ldots,H_{2m-1} \}$ "--- разбиение $Z_2^{m}$ на расширенные $1$"=совершенные коды из определения $\Phi$,
$0\in H_0$, $H_1$ "--- нечетновесовой код,
$H_2$ "--- четновесовой код, $\bar 0\not\in H_2$.
\end{pro}
\begin{proof}
Утверждение следует почти непосредственно из определения отображения $\Phi$.
Действительно, каждое кодовое слово $\bar z=(z_1,\ldots,z_n)\in {\cal C}$ добавляет
$SW_{z_1}\cdot\ldots\cdot SW_{z_n}$ к $SW_{\cal C}(X,Z,T)$, где
$$\begin{array}{rcl@{\qquad}l}
  SW_{0} &\triangleq& X, \\
  SW_{2j+1}  &\triangleq& Z,& j=0, \ldots ,m-1,\\
  SW_{2j}  &\triangleq& T, & j=1, \ldots ,m-1.
\end{array}$$
С другой стороны, как следует из определения отображения $\Phi$,
каждое кодовое слово $\bar z=(z_1,\ldots,z_n)\in {\cal C}$ добавляет
$W_{H_{z_1}}(X,Y)\cdot\ldots\cdot W_{H_{z_n}}(X,Y)$ к $W_{\tilde C}(X,Y)$.
Соотношения
$$\begin{array}{rcl@{\qquad}l}
  W_{H_{2j+1}}(X,Y)&=&W_{H_{1}}(X,Y), & j=0, \ldots ,m-1,\\
  W_{H_{2j}}(X,Y)&=&W_{H_{2}}(X,Y), & j=1, \ldots ,m-1,
\end{array}$$
вытекающие из \hyperref[cl_ru:W_H]{Леммы~}\ref{cl_ru:W_H}, завершают доказательство.
\end{proof}

Следующая теорема "--- основной результат этого раздела.
\begin{theorem}\label{cl_ru:duality}
Пусть ${\cal C}$ и ${\cal C}^\bot$ "--- дуальные линейные $Z_{2^k}$-коды,
$C=\varphi({\cal C})$ и $\tilde C_\bot=\Phi({\cal C}^\bot)$.
Тогда коды $C$ и $\tilde C_\bot$ являются формально дуальными, т.\,е.
\[
 W_C(X,Y)=\frac{1}{|\tilde C_\bot|}W_{\tilde C_\bot}(X+Y,X-Y).
\]
\end{theorem}
\begin{proof}
Согласно \hyperref[cl_ru:CC]{Леммам~}\ref{cl_ru:CC} \hyperref[cl_ru:W_H]{и~}\ref{cl_ru:W_H}
$$\begin{array}{r@{\,}c@{\,}l}
W_C(X,Y)
&=&\frac{1}{|{\cal C}^\bot|}SW_{{\cal C}^\bot} \bigg(
\frac{1}{|H_0|}W_{H_0}(X{+}Y,X{-}Y),
\frac{1}{|H_1|}W_{H_1}(X{+}Y,X{-}Y),
\frac{1}{|H_2|}W_{H_2}(X{+}Y,X{-}Y)
\bigg)\\
&=&\frac{1}{|{\cal C}^\bot|}\left(\frac{2m}{2^m}\right)^{\!\!n}SW_{{\cal C}^\bot} \Big(
W_{H_0}(X{+}Y,X{-}Y),
W_{H_1}(X{+}Y,X{-}Y),
W_{H_2}(X{+}Y,X{-}Y)
\Big).
\end{array}$$
Остается заметить, что по \hyperref[cl_ru:DcoLin]{Лемме~}\ref{cl_ru:DcoLin} мы имеем
${|{\cal C}^\bot|}\left(\frac{2^m}{2m}\right)^n =|\tilde C_\bot|$
и по \hyperref[cl_ru:W-SW]{Предложению~}\ref{cl_ru:W-SW} получаем
$$
SW_{{\cal C}^\bot} \Bigl(W_{H_0}(X{+}Y,X{-}Y),
W_{H_1}(X{+}Y,X{-}Y),
W_{H_2}(X{+}Y,X{-}Y){\Bigr)}
{}={}W_{\tilde C_\bot}(X{+}Y,X{-}Y).
$$
\end{proof}
В \hyperref[cl_ru:S_perf]{Разделах~}\ref{cl_ru:S_perf} \hyperref[cl_ru:S_Had]{и~}\ref{cl_ru:S_Had}
мы построим два класса кодов, $Z_{2^k}$-дуальных друг другу:
ко-$Z_{2^k}$"=линейные расширенные $1$"=совершенные коды и $Z_{2^k}$"=линейные коды Адамара.
В \hyperref[cl_ru:S_nonexist]{последнем разделе} мы докажем полноту конструкции.
\section[Ко-$Z_{2^k}$"=линейные расширенные $1$"=совершенные коды]{%
Ко-$Z_{2^k}$"=линейные расширенные $1$"=совершенные коды}
\label{cl_ru:S_perf}

{Э}{тот} раздел посвящен расширенным $1$"=совершенным кодам.
%
Сначала мы вводим понятие {\em \1"=совершенного кода}, являющееся обобщением понятия расширенного $1$"=совершенного кода на некоторые недвоичные случаи.
Как и для случая $1$"=совершенных кодов, вопрос существования \1"=совершенных кодов в различных пространствах интересен сам по себе.
Затем в \hyperref[cl_ru:ss_b]{Подразделе~}\ref{cl_ru:ss_b}
мы строим класс {\1"=совершенных кодов} в $(Z_{2^k}^n,d^\diamond)$.
В \hyperref[cl_ru:ss_c]{Подразделе~}\ref{cl_ru:ss_c} мы подводим итог:
образы таких кодов при отображении $\Phi$ являются ко-$Z_{2^k}$"=линейными расширенными $1$"=совершенными кодами.
В \hyperref[cl_ru:ss_d]{Подразделе~}\ref{cl_ru:ss_d}
мы приводим примеры \1"=совершенных кодов в $Z_{2m}$, где $m$ не есть степень двойки.

\subsection[$1'$-Совершенные коды]{\1-Совершенные коды} \label{cl_ru:ss_a}
Пусть $G=(V,E)$ "--- регулярный двудольный граф с долями $V_{ev}$, $V_{od}$.
Подмножество $C\subseteq V_{ev}$ называется {\em \1"=совершенным кодом}, если 
для каждого $\bar x\in V_{od}$ найдется ровно один смежный $\bar x$ элемент $\bar c\in C$.
Нетрудно видеть, что \1"=совершенный код является оптимальным кодом с расстоянием $4$,
т.\,е. его мощность максимальна по всем кодам с расстоянием $4$ в том же пространстве.
В частности, такой код является $3$-диаметрально совершенным кодом в смысле \cite{cl_ru:AhlAydKha}.
\1-Совершенные коды в двоичном пространстве Хемминга известны как расширенные $1$"=совершенные коды.
Следующий критерий можно взять в качестве альтернативного определения  \1"=совершенного кода.
\begin{pro}\label{cl_ru:1perf}
Пусть $G=(V,E)$ "--- регулярный двудольный граф степени $r>0$
и $d^{G}$ "--- естественная графская метрика на $V$.
Множество вершин $C\subset V$ является \1"=совершенным кодом тогда и только тогда, когда
$|C|=|V|/2r$ и $d^{G}(\bar c_1,\bar c_2)\geq 4$
для любых различных $\bar c_1, \bar c_2\in C$.
\end{pro}
\begin{proof}
{\em Только если:} Пусть $C\subset V$ есть \1"=совершенный код.
По определению он является подмножеством некоторой доли $V_{ev}$
графа $(V,E)$ и расстояние между любыми двумя элементами кода $C$ четно.

С другой стороны, если $d^{G}(\bar c_1, \bar c_2)=2$ для некоторых $\bar c_1, \bar c_2\in C$,
то найдется такой элемент $\bar x\in V$, что
$d^{G}(\bar x, \bar c_1)=1$ и $d^{G}(\bar x, \bar c_2)=1$.
Существование такого элемента противоречит определению \1"=совершенного кода.
Следовательно, $d^{G}(\bar c_1, \bar c_2)\geq 4$ для любых различных $\bar c_1, \bar c_2\in C$.

Каждый элемент множества $V_{od}\triangleq V\backslash V_{ev}$ является смежным ровно с одним элементом кода $C$.
С другой стороны, каждый элемент $C$ смежный ровно с $r$ элементами $V_{od}$.
Отсюда $|C|=|V_{od}|/r=|V|/2r$.

{\em Если:}
Предположим, что $C$ имеет мощность $|V|/2r$ и кодовое расстояние не меньше $4$.
Обозначим
$$
C_1 \triangleq \{\bar x\,|\,d^{G}(\bar x,C)=1\}.
$$
Каждый элемент из $C_1$ имеет ровно одного соседа из $C$
(в противном случае кодовое расстояние $C$ не превышает $2$).
Следовательно, $|C_1|=r|C|=|V|/2$.

(*) {\em Мы утверждаем, что $C_1$ не содержит двух смежных элементов.}
Допустим противное.
Тогда граф $G$ содержит цепь $(\bar c, \bar x, \bar y, \bar c')$,
где $\bar c, \bar c'\in C$, $\bar x, \bar y \in C_1$.
Случай $\bar c=\bar c'$ противоречит двудольности графа $G$.
Случай $\bar c\neq\bar c'$ противоречит кодовому расстоянию кода $C$.
Утверждение (*) доказано.

Поскольку $G$ "--- регулярный граф, другая половина вершин $V \backslash C_1$ также не содержит смежных элементов.
Следовательно, $C$ является \1"=совершенным кодом по определению, где $V_{od}=C_1$.
\end{proof}

\subsection[Класс $1'$"=совершенных кодов в $Z_{2^k}^n$]{Класс \1"=совершенных кодов в $Z_{2^k}^n$} \label{cl_ru:ss_b}
Пусть $n=2^r$ и $I=(i_1,\ldots,i_k)$ "--- набор неотрицательных целых чисел, удовлетворяющий равенству
$1 i_1+2 i_2+ \ldots +k i_k=r$.
Пусть $\bar b_1,\ldots,\bar b_n\in Z_{2^k}^{1+i_1+...+i_k}$ "--- все элементы множества
 $\{1\}\times (2^{k-1} Z_{2^k})^{i_1}\times (2^{k-2} Z_{2^k})^{i_2}
\times \ldots \times (2^{0} Z_{2^k})^{i_k}$, упорядоченные лексикографически.
И пусть $B_I$ "--- матрица со столбцами $\bar b_1,\ldots,\bar b_n$.

\begin{examp}
При $k=3$, $i_1=2$, $i_2=1$, $i_3=0$ имеем $r=4$, $n=2^r=16$ и
$$
B_I=
\left(%
\begin{array}{c}
1111111111111111 \\
0000000044444444 \\
0000444400004444 \\
0246024602460246 \\
\end{array}
\right).
$$
При $k=3$, $i_1=0$, $i_2=i_3=1$ имеем $r=5$, $n=2^r=32$ и
$$
B_I=
\left(%
\begin{array}{c}
11111111111111111111111111111111 \\
00000000222222224444444466666666 \\
01234567012345670123456701234567 \\
\end{array}
\right).
$$
При $k=3$, $i_1=i_2=i_3=1$ имеем $r=6$, $n=2^r=64$ и
$$
B_I=
\left(%
\begin{array}{c}
1111111111111111111111111111111111111111111111111111111111111111 \\
0000000000000000000000000000000044444444444444444444444444444444 \\
0000000022222222444444446666666600000000222222224444444466666666 \\
0123456701234567012345670123456701234567012345670123456701234567 \\
\end{array}
\right).
$$
\end{examp}

\begin{lemma}\label{cl_ru:isPerf}
Линейный код
${\cal H}_I\triangleq \{\bar h=(h_1,\ldots,h_n)\in Z_{2^k}^n\,
|\,\sum_{j=1}^n h_j \bar b_j = B_I\bar h^T = \bar 0 \}$ с проверочной матрицей $B_I$ является \1"=совершенным.
\end{lemma}
\begin{proof}
Сначала покажем, что

(*) {\em расстояние $2$ между кодовыми словами
$\bar h_1$, $\bar h_2$ из ${\cal H}_I$ невозможно}.
Действительно, если $d^\diamond(\bar h_1,\bar h_2)=2$, то кодовое слово
$\bar h \triangleq \bar h_1-\bar h_2$ имеет одну или две ненулевых позиции.
Первый случай противоречит факту $\bar b_i\neq \bar 0$, $i=1, \ldots ,n$.
Во вором случае мы имеем $\beta \bar b_i+\gamma \bar b_j=\bar 0$ для некоторых различных $i$, $j$
и нечетных $\beta$, $\gamma$. Но первый ряд матрицы $B_I$ подразумевает, что $\beta=-\gamma$ и,
поскольку $\beta$ и $\gamma$ нечетные, мы получаем $\bar b_i=\bar b_j$.
Это опять противоречит построению матрицы $B_I$. Утверждение (*) доказано.

Первый ряд матрицы $B_I$ подразумевает, что слова ${\cal H}_I$ имеют четный вес.

Нам нужно проверить, что каждое нечетное слово находится на расстоянии $1$ ровно от одного кодового слова.
Пусть $\bar z=(z_1, \ldots ,z_n)\in Z_{2^k}^n$ имеет нечетный вес,
$\bar s\triangleq\sum_{j=1}^n z_i\bar b_i$ и $s\triangleq\sum_{j=1}^n z_i$.
Заметим, что $s$ нечетно.
Поскольку $s^{-1}\bar s\in \{1\}\times (2^{k-1} Z_{2^k})^{i_1}\times (2^{k-2} Z_{2^k})^{i_2}
\times \ldots \times (2^{0} Z_{2^k})^{i_k}$, найдется такое $j'$, что
$s^{-1}\bar s=\bar b_{j'}$. Пусть $\bar z'\in Z_{2^k}^n$ "--- слово с $s$ в $j'$й позиции и нулями в остальных.
Легко проверить, что $\bar z-\bar z'$ "--- кодовое слово на расстоянии $1$ от $\bar z$.
Как следует из (*), такое кодовое слово единственно для каждого нечетного $\bar z$.
\end{proof}

\subsection{Ко-$Z_{2^k}$"=линейные расширенные $1$"=совершенные двоичные коды} \label{cl_ru:ss_c}

Теперь у нас есть все необходимое для построения класса ко-$Z_{2^k}$"=линейных расширенных $1$"=совершенных кодов.
Как мы увидим в \hyperref[cl_ru:S_nonexist]{Разделе~}\ref{cl_ru:S_nonexist},
построенный класс исчерпывает все такие коды при условии, что отображение $\Phi$ фиксировано.

\begin{theorem}\label{cl_ru:perf}
Ко-$Z_{2^k}$"=линейный код $\tilde H_I\triangleq \Phi({\cal H}_I)$ является двоичным\linebreak[5] $(nm,2^{nm}/2nm,4)$-кодом, т.\,е. расширенным $1$"=совершенным кодом.
\end{theorem}

\begin{proof}
Утверждение следует непосредственно из
\hyperref[cl_ru:DcoLin]{Лем\-мы~}\ref{cl_ru:DcoLin},
\hyperref[cl_ru:DcoLin]{Пред\-ло\-же\-ния~}\ref{cl_ru:1perf} и
\hyperref[cl_ru:DcoLin]{Лем\-мы~}\ref{cl_ru:isPerf}.
\end{proof}

\subsection{Замечание об общем случае $Z_{2m}$} \label{cl_ru:ss_d}

На самом деле, \1"=совершенные коды могут быть построены над $Z_{2m}$ с расстоянием $d^\diamond$ для каждого $m\geq 1$, см. приведенные ниже примеры.
Поскольку $m=2^\mu$ есть необходимое условие для построения двоичных кодов способом, описанным в \hyperref[cl_ru:S_co-Z-lin]{Подразделе~}\ref{cl_ru:S_co-Z-lin},
классификация \1"=совершенных кодов в других случаях выходит за рамки настоящей работы.

\begin{examp}\label{cl_ru:ex1perf} Коды с проверочными матрицами

$$B'\triangleq
\left(%
\begin{array}{cccccccc}
  1 & 1 & 1 & 1 & 1 & 1 & 1 & 1  \\
  0 & 0 & 0 & 0 & 12& 12& 12& 12 \\
  0 & 6 & 12& 18& 0 & 6 & 12& 18
\end{array}%
\right),
$$
$$
B''\triangleq
\left(%
\begin{array}{cccccc}
  1 & 1 & 1  & 1  & 1  & 1  \\
  0 & 4 & 8  & 12 & 16 & 20
\end{array}%
\right)
$$
являются \1"=совершенными кодами над $Z_{24}$ с расстоянием $d^\diamond$.
\end{examp}

\section[$Z_{2^k}$-Линейные коды Адамара]{%
$Z_{2^k}$-Линейные коды Адамара}\label{cl_ru:S_Had}

\begin{lemma} \label{cl_ru:L_HD}
Линейный код ${\cal H}$ в $Z_{2^k}^n$
является $(n,n2^{k},n2^{k-2})^*$-кодом
тогда и только тогда, когда ${\cal H}^\bot$ "---
 \1"=совершенный код в $(Z_{2^k}^n,d^\diamond)$.
\end{lemma}
\begin{proof}
Предположим, что линейный код ${\cal H}$ имеет параметры $(n,n2^{k},n2^{k-2})^*$.
Тогда $\varphi({\cal H})$ имеет параметры кода Адамара, см. \hyperref[cl_ru:Cor1]{Следствие~}\ref{cl_ru:Cor1}. Согласно \hyperref[cl_ru:duality]{Теоремы~}\ref{cl_ru:duality}
код $\Phi({\cal H}^\bot)$ является формально дуальным коду Адамара $\varphi({\cal H})$, т.\,е.
$\Phi({\cal H}^\bot)$ "--- расширенный $1$"=совершенный код.
Тогда мощность и минимальное $d^\diamond$-расстояние $4$ кода ${\cal H}^\bot$ следуют из \hyperref[cl_ru:DcoLin]{Леммы~}\ref{cl_ru:DcoLin}.

Обратное утверждение доказывается обратным рассуждением.
\end{proof}

Следующая теорема следует сразу из \hyperref[cl_ru:isPerf]{Лемм~}\ref{cl_ru:isPerf}
\hyperref[cl_ru:L_HD]{и~}\ref{cl_ru:L_HD}.

\begin{theorem}
Код
${\cal D}_I\triangleq {\cal H}_I^\bot$
"--- линейный $(n,n2^{k},n2^{k-2})^*$-код с порождающей матрицей $B_I$.
Соответствующий $Z_{2^k}$-линейный код $D_I\triangleq \varphi({\cal D}_I)$ является двоичным $(n2^{k-1},n2^{k},n2^{k-2})$-кодом, т.\,е. кодом Адамара.
\end{theorem}

\begin{note}
Код ${\cal D}_{(r,0,\ldots,0)}$ известен как
{\em код Рида-Маллера первого порядка над $Z_{2^k}$} \cite{cl_ru:GuptaOth2005}.
\end{note}

\begin{note} Если существует код Адамара $A$ длины $m$ (см. \hyperref[cl_ru:S_Z-lin]{Подраздел~}\ref{cl_ru:S_Z-lin}),
то $Z_{2m}$"=линейный код Адамара длины $nm=2^r m$ может быть построен даже если $m$ не является степенью двойки.
Например, линейный $Z_{24}$-код с порождающей матрицей $B'$
из \hyperref[cl_ru:ex1perf]{Примера~}\ref{cl_ru:ex1perf}
имеет параметры $(16,192,48)^*$ и соответствующий двоичный $Z_{24}$"=линейный код имеет параметры $(96,192,48)$, т.\,е. параметры кода Адамара длины $96$.

Однако код с порождающей матрицей $B''$ (\hyperref[cl_ru:ex1perf]{Пример~}\ref{cl_ru:ex1perf})
имеет параметры $(6,144,30)^*$, соответствующие двоичному $(72,144,30)$-коду.
Кодовое расстояние этого кода меньше, чем кодовое расстояние кода Адамара той же длины и той же мощности.
\end{note}

\section[Несуществование неизвестных ко-$Z_{2^k}$"=линейных расширенных совершенных кодов и $Z_{2^k}$"=линейных кодов Адамара]{%
Несуществование неизвестных ко-$Z_{2^k}$"=линейных расширенных совершенных кодов и $Z_{2^k}$"=линейных кодов Адамара}\label{cl_ru:S_nonexist}
{П}{усть} $n$ "--- степень двойки. В этом разделе мы покажем (\hyperref[cl_ru:Th_HH]{Теорема~}\ref{cl_ru:Th_HH}), что каждый линейный $(n, 2^{kn}/n2^k, 4)^\diamond$-код в $Z_{2^k}^n$ эквивалентен некоторому коду из класса, построенного в \hyperref[cl_ru:S_perf]{Разделе~}\ref{cl_ru:S_perf}. Аналогично, каждый линейный
$(n,n2^{k},n2^{k-2})^*$-код в $Z_{2^k}^n$ эквивалентен некоторому коду из класса, построенного в \hyperref[cl_ru:S_Had]{Разделе~}\ref{cl_ru:S_Had}
(\hyperref[cl_ru:Th_DD]{Теорема~}\ref{cl_ru:Th_DD}). Ключевым моментом доказательства является \hyperref[cl_ru:Th_HH]{Лемма~}\ref{cl_ru:L_1111}.
Частный $Z_4$-случай этой леммы был доказан в \cite{cl_ru:Kro:2000:Z4_Perf} и в \cite{cl_ru:BorPheRif:2003}, приведенное ниже доказательство использует другой ход рассуждений.

Мы говорим, что два линейных кода ${\cal C}_1,{\cal C}_2\subseteq Z_{2^k}^n$  \emph{эквивалентны}, если 
${\cal C}_2=\bar z\circ\pi {\cal C}_1$,
где $\pi$ "--- перестановка координат,
$\bar z \in (Z_{2^k}^*)^n\triangleq \{1,3,\ldots,2^k-1\}^n$,
и $\circ$ "--- покоординатное произведение, т.\,е.
$(z_1,\ldots,z_n)\circ(x_1,\ldots,x_n)\triangleq (z_1 x_1,\ldots,z_n x_n)$.
Заметим, что как $\pi$, так и $\bar z\,\circ$ "--- групповые автоморфизмы аддитивной группы $Z_{2^k}^n$ и изометрии метрических пространств
$(Z_{2^k}^n,d^*)$ и $(Z_{2^k}^n,d^\diamond)$.
Следующее предложение показывает, что нет других линейных изометрий пространства $(Z_{2^k}^n,d^*)$ или $(Z_{2^k}^n,d^\diamond)$, следовательно наше определение эквивалентности естественно.
\begin{pro}
Предположим, что $\Gamma$ является линейным преобразованием $Z_{2^k}^n$ и изометрией $(Z_{2^k}^n,d)$, где $d=d^*$ или $d=d^\diamond$.
Тогда $\Gamma(\bar x)\equiv\bar z\circ\pi(\bar x)$,
где $\pi$ "--- перестановка координат и $\bar z \in (Z_{2^k}^*)^n$.
\end{pro}
\begin{proof}
Обозначим через $\bar e_i$ слово с единицей в позиции с номером $i$ и нулями в остальных.
Достаточно проверить, что для каждого $i$ мы имеем
$\Gamma(\bar e_i)\equiv z_j \bar e_{j}$ с некоторыми
$j=\pi(i)$ и $z_{j}\in Z_{2^k}^*$.
Действительно, из изометрических свойств $\Gamma$
мы выводим, что $\Gamma(\bar e_i)$ имеет только одну ненулевую координату.
С другой стороны, значение этой координаты принадлежит $Z_{2^k}^*$,
поскольку $\Gamma$ является биекцией.
\end{proof}

Следующие две теоремы являются основным результатом этого раздела.
Напомним, что коды ${\cal H}_I$ и ${\cal D}_I$ определяются в
\hyperref[cl_ru:S_perf]{Разделах~}\ref{cl_ru:S_perf} \hyperref[cl_ru:S_Had]{и~}\ref{cl_ru:S_Had}.

\begin{theorem} \label{cl_ru:Th_HH}
Пусть ${\cal H}\subset Z_{2^k}^n$ "--- некоторый линейный
\1"=совершенный код.
Тогда ${\cal H}$ эквивалентен ${\cal H}_I$ для некоторого набора
$I=(i_1,\ldots,i_k)$ неотрицательных целых чисел, удовлетворяющего равенству $1 i_1+2 i_2+ \ldots +k i_k=\log_2 n$.
\end{theorem}

\begin{theorem} \label{cl_ru:Th_DD}
Пусть ${\cal D}\subset Z_{2^k}^n$ "--- линейный $(n,n2^{k},n2^{k-2})^*$-код.
Тогда ${\cal D}$ эквивалентен ${\cal D}_I$ для некоторого набора
$I=(i_1,\ldots,i_k)$ неотрицательных целых чисел, удовлетворяющего равенству
$1 i_1+2 i_2+ \ldots +k i_k=\log_2 n$.
\end{theorem}

Согласно \hyperref[cl_ru:L_HD]{Леммы~}\ref{cl_ru:L_HD}
достаточно доказать только одну из \hyperref[cl_ru:Th_HH]{Теорем~}\ref{cl_ru:Th_HH}, \ref{cl_ru:Th_DD}.
Нам понадобится следующий вспомогательный результат.

\begin{lemma}\label{cl_ru:L_1111}
Если $\cal D$ "--- линейный $(n,n2^{k},n2^{k-2})^*$-код в $Z_{2^k}^n$, то
$\cal D$ содержит элемент из $(Z_{2^k}^*)^n$.
\end{lemma}
\begin{proof}
Докажем лемму по индукции.
Если $n=1$,
то ${\cal D}=Z_{2^k}$.
Предположим, что $n>1$.

(*) \emph{Мы утверждаем, что $\cal D$ содержит элемент из $\{0,2^{k-1}\}^n$ веса $n2^{k-2}$.}
Пусть $\bar x'$ и $\bar x''$ "--- два линейно независимых элемента в $\cal D$ порядков $2m'$ и $2m''$ соответственно.
Это означает, что $m' \bar x'$ и $m'' \bar x''$ "--- различные ненулевые элементы из ${\cal D} \cap \{0,2^{k-1}\}^n$.
Поскольку $\varphi$ является изометричным вложением (см. \hyperref[cl_ru:P_iso]{Предложение~}\ref{cl_ru:P_iso})
и $\varphi({\cal D})$ есть $(n2^{k-1},n2^{k},n2^{k-2})$-код Адамара,
единственные возможные значения $wt^*$-весов элементов $\cal D$ есть $0$, $n2^{k-2}$ и $n2^{k-1}$.
Таким образом, как минимум один из $m' \bar x'$ и $m'' \bar x''$ имеет вес $n2^{k-2}$.
Утверждение (*) доказано.

Без потери общности будем считать, что
\(\bar a \triangleq (2^{k-1},\ldots,2^{k-1},0,\ldots,0)
\in {\cal D}.\)

Рассмотрим два кода, полученных из $\cal D$ выкалыванием $n/2$ координат:
$$\begin{array}{rcl}
{\cal D}_1 & {}\triangleq{} & \{\bar z'\in Z_{2^k}^{n/2} \, | \  %
\exists\, \bar z''\in Z_{2^k}^{n/2} :\, (\bar z',\bar z'')\in {\cal D}\}, \\
{\cal D}_2 & {}\triangleq{} & \{\bar z''\in Z_{2^k}^{n/2} \, | \  %
\exists\, \bar z'\in Z_{2^k}^{n/2} :\, (\bar z',\bar z'')\in {\cal D}\}.
\end{array}$$

(**) \emph{Мы утверждаем, что ${\cal D}_1$ и ${\cal D}_2$ "--- линейные $(n/2,n2^{k-1},n2^{k-3})^*$-коды.}
Линейность очевидна.
Кодовое расстояние следует из факта, что расстояние между $\bar a$
и любым элементом $\bar z\in {\cal D}$ равно $0$, $n2^{k-2}$ или $n2^{k-1}$.
Мощность кода ${\cal D}_1$ удовлетворяет неравенству $|{\cal D}_1| \geq |{\cal D}|/2$, поскольку кодовое расстояние ${\cal D}$ допускает только одно ненулевое слово с нулями в первых $n/2$ координатах.
С другой стороны, есть такое кодовое слово:
$(0,\ldots,0,2^{k-1},\ldots,2^{k-1})=(2^{k-1},\ldots,2^{k-1})+\bar a\in {\cal D}$.
Поэтому $|{\cal D}_1| = |{\cal D}|/2 = n2^{k-1}$.
Аналогично $|{\cal D}_2| = n2^{k-1}$. Утверждение (**) доказано.

По предположению индукции ${\cal D}_1$ содержит элемент $\bar u'$
из $(Z_{2^k}^*)^{n/2}$.
Это означает, что существует $\bar u''$ из $Z_{2^k}^{n/2}$ такое, что
$(\bar u',\bar u'')\in {\cal D}$.
Поскольку $wt^*(2^{k-1}(\bar u',\bar u''))\in \{n2^{k-2},n2^{k-1}\}$ и
$2^{k-1}\bar u'=(2^{k-1},...,2^{k-1})$,
имеем $wt^*(2^{k-1}\bar u'')=0$ или $wt^*(2^{k-1}\bar u'')=n2^{k-2}$.
Таким образом,
\setlength{\arraycolsep}{0mm}
\begin{eqnarray}
 \label{cl_ru:u''0} 2^{k-1}\bar u'' & {}={} & (0,\ldots,0), \\
 \label{cl_ru:u''1} \mbox{или}\ \ 2^{k-1}\bar u'' & {}={} & (2^{k-1},\ldots,2^{k-1}).
\end{eqnarray}
Аналогично, рассматривая код ${\cal D}_2$, мы можем найти кодовое слово $(\bar v',\bar v'')\in {\cal D}$ такое, что
$2^{k-1}\bar v''=(2^{k-1},...,2^{k-1})$ и $v'$ удовлетворяет
\setlength{\arraycolsep}{0mm}
\begin{eqnarray}
 \label{cl_ru:v'0} 2^{k-1}\bar v' & {}={} & (0,\ldots,0), \\
 \label{cl_ru:v'1} \mbox{или}\ \ 2^{k-1}\bar v' & {}={} & (2^{k-1},\ldots,2^{k-1}).
\end{eqnarray}
Если выполнено (\ref{cl_ru:u''1}), то $(u',u'')\in (Z_{2^k}^*)^{n}\cap {\cal D}$.
Если выполнено (\ref{cl_ru:v'1}), то $(v',v'')\in (Z_{2^k}^*)^{n}\cap {\cal D}$.
Если выполнено (\ref{cl_ru:u''0}) и (\ref{cl_ru:v'0}), то $(u',u'')+(v',v'')\in (Z_{2^k}^*)^{n}\cap {\cal D}$.
\hyperref[cl_ru:L_1111]{Лемма~}\ref{cl_ru:L_1111} доказана.
\end{proof}

\begin{proof}[{\hyperref[cl_ru:Th_DD]{Теоремы~}\ref{cl_ru:Th_DD}}]
По \hyperref[cl_ru:L_1111]{Лемме~}\ref{cl_ru:L_1111} код
$\cal D$ содержит элемент $\bar c=(c_1,\ldots,c_n)$ из $(Z_{2^k}^*)^n$.
Тогда код ${\cal D}'\triangleq (c_1^{-1},\ldots,c_n^{-1})\circ {\cal D}$
эквивалентен коду ${\cal D}$ и содержит $\bar 1=(1,\ldots,1)$.
Пусть $\{\bar 1, q_1,\ldots, q_s\}$ "--- некоторый базис кода ${\cal D}'$ и
$i_j$ "--- число элементов порядка $2^{j}$, $j=1,\ldots,k$.
Пусть $\bar 1, q_1,\ldots, q_s$ "--- строки матриц $Q$
размера $(1+i_1+\ldots+i_k)\times n$, порождающей матрицы кода ${\cal D}'$.
Тогда столбцы $Q$ являются элементами $\{1\}\times (2^{k-1} Z_{2^k})^{i_1}\times (2^{k-2} Z_{2^k})^{i_2}
\times \ldots \times (2^{0} Z_{2^k})^{i_k}$.
Поскольку по \hyperref[cl_ru:L_HD]{Лемме~}\ref{cl_ru:L_HD} кодовое $d^\diamond$-расстояние кода ${\cal D}'^\bot$ больше $2$,
все столбцы попарно различны.
Следовательно $Q$ совпадает с $B_I$, $I=(i_1,\ldots,i_k)$,
с точностью до перестановки столбцов.
\end{proof}

Таким образом, мы заключаем, что при условии фиксированных отображений $\varphi$ и $\Phi$ все с точностью до эквивалентности
ко-$Z_{2^k}$"=линейные расширенные $1$"=совершенные коды и
$Z_{2^k}$"=линейные коды Адамара описаны в \hyperref[cl_ru:S_perf]{Разделах~}\ref{cl_ru:S_perf} \hyperref[cl_ru:S_Had]{и~}\ref{cl_ru:S_Had}.

\ifx\href\undefined \newcommand\href[2]{#2} \fi\ifx\url\undefined
  \newcommand\url[1]{\href{#1}{#1}} \fi